\newif\ifincludeauthorship
\includeauthorshiptrue

\newif\ifrefereeview
\refereeviewfalse

\ifrefereeview
  \documentclass[12pt]{amsart}
\else
  \documentclass{amsart}
\fi

\usepackage{amsfonts,amsmath,amsthm,amssymb,caption,color,enumitem,float,graphicx,nicefrac,subcaption,url,wrapfig}
\usepackage{fullpage,lineno,setspace}

\usepackage[toc,page]{appendix}
\usepackage[noadjust]{cite}
\usepackage[foot]{amsaddr}
\usepackage[mathscr]{eucal} \usepackage{mathrsfs}
\usepackage{tikz}
\usetikzlibrary{automata,arrows,positioning,calc}
\usepackage[mathscr]{eucal}

\DeclareGraphicsExtensions{.png,.jpeg,.jpg,.pdf}
\graphicspath{{./Figures/}}
\allowdisplaybreaks

\theoremstyle{plain}
\newtheorem{theorem}{Theorem}[section]

\newtheorem{proposition}[theorem]{Proposition}
\newtheorem{lemma}[theorem]{Lemma}

\newtheorem{remark}[theorem]{Remark}

\newcommand{\lb}{\left\{}
\newcommand{\rb}{\right\}}
\newcommand{\Def}{\overset{\textbf{def}}{=}}
\newcommand{\R}{\mathbb{R}}
\newcommand{\N}{\mathbb{N}}
\newcommand{\Q}{\mathbb{Q}}
\newcommand{\eps}{\varepsilon}
\newcommand{\bOne}{\mathbf{1}}
\newcommand{\filt}{\mathscr{F}}
\newcommand{\BP}{\mathbb{P}}
\newcommand{\BE}{\mathbb{E}}
\newcommand{\uphi}{\underline{\phi}}
\newcommand{\KK}{\textbf{K}}
\newcommand{\be}{\mathbf{e}}
\newcommand{\Space}{\mathbf{S}}
\newcommand{\Region}{\mathsf{R}}
\newcommand{\Upper}{\mathsf{U}}
\newcommand{\Lower}{\mathsf{L}}
\newcommand{\BB}{\mathbf{B}}
\newcommand{\la}{\left\langle}
\newcommand{\ra}{\right\rangle}
\newcommand{\bdy}{\partial}
\newcommand{\MID}{\mathcal{T}}
\newcommand{\ff}{\mathbf{f}}
\newcommand{\yy}{\mathfrak{y}}
\newcommand{\Agree}{\mathsf{A}}
\newcommand{\agree}{\mathsf{a}}

\newcommand{\smooth}{\mathsf{s}}
\newcommand{\vsig}{\varsigma}

\author{Hossein Nick Zinat Matin}
\address{Department of Civil and Environmental Engineering,\\
    University of California at Berkeley\\
    Berkeley, CA 94720}
\email{h-matin@berkeley.edu (nickzin2@illinois.edu)}

\author{Richard B. Sowers}
\address{Department of Industrial and Enterprise Systems Engineering \\
Department of Mathematics\\
    University of Illinois at Urbana--Champaign\\
    Urbana, IL 61801}
\email{r-sowers@illinois.edu}
\thanks{This material is based upon work supported by the National Science Foundation under Grant Number  CMMI 1727785. Part of this research was performed while the author was visiting  the Institute for Pure and Applied mathematics (IPAM), which is supported by the National Science Foundation (Grant No. DMS-1925919).  This research was also supported by the Campus Research Board of the University of Illinois at Urbana-Champaign.  This work is based on the Ph.D. thesis of H.M. }
\date{\today}

\title{Near-collision dynamics in a noisy car-following model}
\begin{document}
\maketitle

\begin{abstract}
We consider a small stochastic perturbation of an optimal velocity car-following model. We give a detailed analysis of behavior near a collision singularity. We show that collision is impossible in a simplified model without noise, and then show that  collision is asymptotically unlikely over large time intervals in the presence of small noise, with the large time interval scaling like the square of the reciprocal of the strength of the noise.  Our calculations depend on careful boundary-layer analyses.
\end{abstract}


\section{Introduction}
The imminent revolution in autonomous driving technologies has led to a wide range of challenges in modelling and understanding the effect of autonomous vehicles.  Air quality and environment \cite{le2021air,stern2019quantifying}, control of bulk traffic \cite{garavello2020multiscale}, fuel usage \cite{chen2019quantifying}, parking \cite{bischoff2018autonomous}, traffic safety \cite{ye2019evaluating}, and a range of other issues \cite{greenblatt2015automated} are likely to be affected by adoption of autonomous vehicles (cf. \cite{garavello2006traffic}).

Our interest here is a mathematical analysis of \emph{stochasticity} in \emph{adaptive cruise control}.  \emph{Platooning} in autonomous vehicles \cite{lu2017platooning} can have significant effect on long-distance vehicle efficiency by managing air resistance and drag \cite{mitra2007pollution}.  Platooning depends on a sequence of vehicles establishing optimal following distances and velocities, and requires establishing tight inter-vehicle controls \cite{kavathekar2011vehicle, zeng2019joint}.

While platooning and car-following are (usually) designed to implement deterministic car-following controls, stochasticity is intrinsic to real physical systems.  Different vehicles in a platoon may have slightly different characteristics, and one might think of randomness across vehicle \lq\lq types".  Secondly, wind and road conditions may generate forces which are best modelled as random processes.  This latter type of stochasticity will be our focus here.

Letting $x$ be the position of a \emph{following} vehicle, which is following a  (rightward-travelling) \emph{preceding} vehicle whose position is $x^{(p)}$, one might build upon currently-accepted models for adaptive cruise control and add noise, with the dynamics (see Section \ref{S:ModelsAndContributions}) as
\begin{equation}
\label{E:OVM}
\begin{aligned}
\ddot x_t & = -\alpha \lb  \dot x_t -  V \left(\frac{x^{(p)}_t - x_t}{d}\right) \rb  - \beta \frac{\dot x_t - \dot x^{(p)}_t}{\left(x^{(p)}_t -x_t\right)^2} + \eps \dot W_t 
\end{aligned}\end{equation} 
where $\alpha$, $\beta$, and $d$ are positive constants, $V$ is an appropriate function (which will be specified in \eqref{E:function_V}), $W$ is a standard Brownian motion, and $\eps$ is a parameter modelling the strength of random forces.  Informally, the first two terms reflect control forces which cause the following vehicle to follow the preceding vehicle at a preferred distance (the first term, with $\alpha$) and to match the velocity $\dot x^{(p)}$ of the preceding vehicle (the second term, with $\beta$).  The preceding vehicle car could in fact be the \emph{lead} car in a \emph{platoon}, or in fact one in a sequence of platooning vehicles (which would correspond to a recursive sequence of equations of the type \eqref{E:OVM}, replacing $x$ with the position $x^{(n)}$ of the $n$-th vehicle in the platoon and replacing $x^{(p)}$ with the position $x^{(n-1)}$ of the $n-1$-st vehicle in the platoon), and perhaps also indexing the Brownian motion by $n$.

Mathematically, the deterministic part of the dynamics of \eqref{E:OVM} reflects a collection of effects; see Section \ref{S:ModelsAndContributions}.   We are interested in the interactions between these effects and noise. Random perturbations can sooner or later drive dynamical systems into a range of states \cite{MR1652127}, and in particular might lead to collisions in a platoon.  We are interested in understanding how this can happen.  We are in particular interested in \emph{small-noise} asymptotics, in the regime where $\eps\searrow 0$ in \eqref{E:OVM}.  Random effects of road and wind are in most cases likely to be small.  Mathematically, small-noise asymptotics often give an understanding of features which are robust to diffusive perturbations \cite{MR1652127}.  

\section{Background and Focus}\label{S:ModelsAndContributions}
\begin{figure}
    \centering
    \includegraphics[width=0.7\columnwidth]{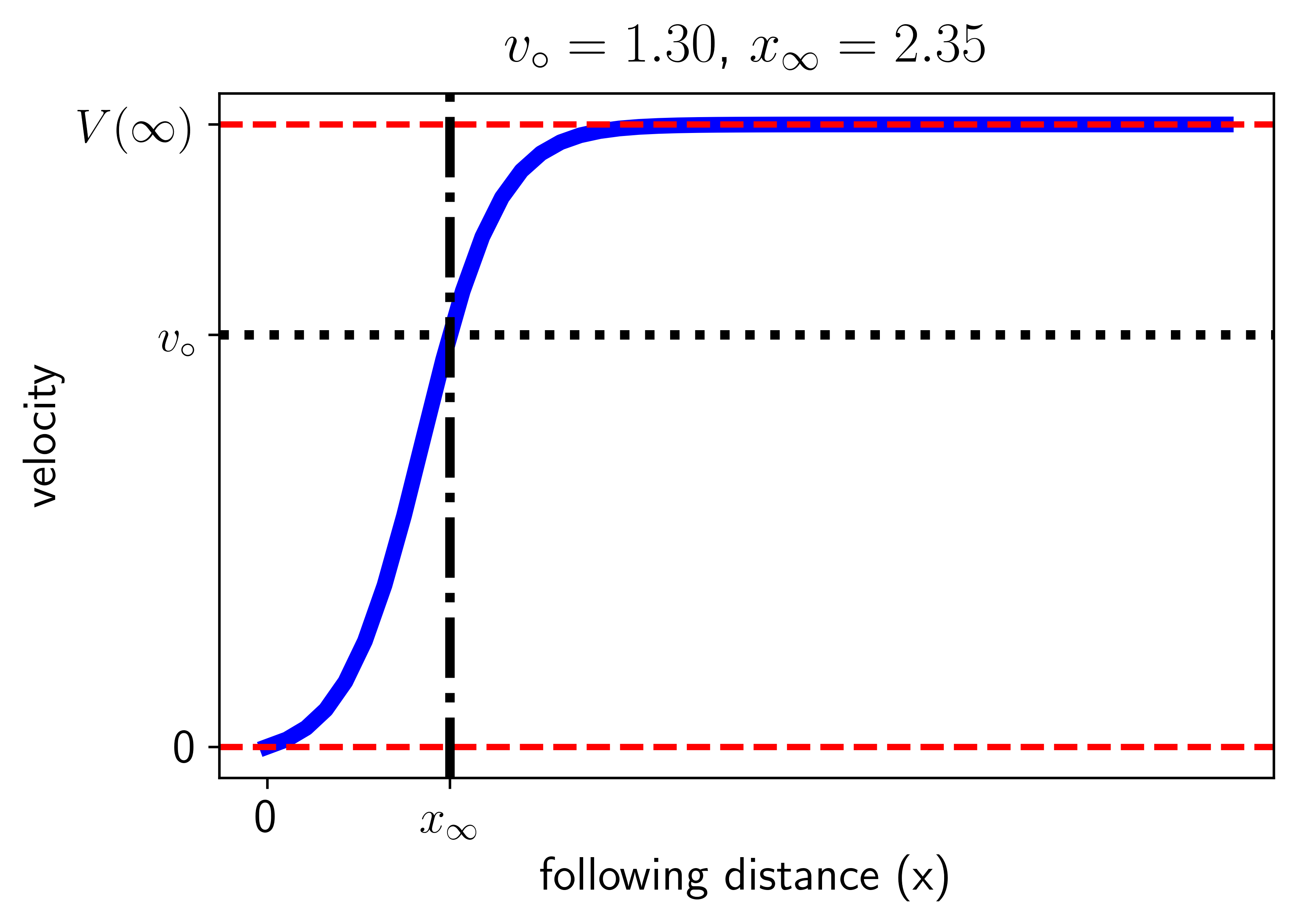}
    \caption{Function $V$}
    \label{fig:Vplot}
\end{figure}

The classic \textit{optimal velocity model}  \cite{bando1994structure, bando1995dynamical}, corresponding to \eqref{E:OVM} with $\beta=\eps=0$, captures the forces needed to follow a preceding vehicle at a fixed distance.  The function $V$, usually required to be smooth, monotone increasing, and bounded, is often taken to be
\begin{equation}\label{E:function_V}
    V(x) \Def \tanh(x-2) - \tanh(-2); \qquad x\in \R 
\end{equation}
see Figure \ref{fig:Vplot}.  The function $V$
transforms the \emph{following distance} $x^{(p)}_t - x_t$ (in a rightward-travelling platoon, $x^{(p)}_t>x_t$) into reference velocity which the following vehicle tries to match. By inverting this transformation (see \eqref{E:equil}), we can get a following distance which is a fixed point for the velocity of the preceding vehicle (and which would represent a fixed point of the dynamics of a platoon). The parameter $d$ captures the sensitivity of the reference velocity with respect to the following distance. The (positive) coefficient $\alpha$ implies that this effect has relaxation time $1/\alpha$.  The $\beta$ term in \eqref{E:OVM}, which alternately has form
\begin{equation*} -\frac{\beta}{\left(x^{(p)}_t -x_t\right)^2}\lb \dot x_t - \dot x^{(p)}_t\rb, \end{equation*}
models a force which directly tries to match the velocity of the preceding vehicle.  This occurs with instantaneous relaxation time $\left( x^{(p)}_t-x_t\right)^2/\beta$; as the following distance becomes small, this relaxation time decreases, and a \emph{singularity} occurs at collision \cite{gazis1961nonlinear}.  The deterministic optimal-velocity follow-the-leader model; i.e., \eqref{E:OVM} with $\eps=0$ combines both of these effects \cite{chiarello2021statistical}.  If $V$ were to be linear, \eqref{E:OVM} would be similar to a proportional-integral-derivative controller \cite{willis1999proportional}.  The $\beta$ term makes the system more robust as it directly models relaxation to the velocity of the preceding vehicle \cite{ helbing1998generalized, jiang2001full, treiber2013traffic}.%

A platoon of vehicles following \eqref{E:OVM} (with $x^{(n)}$ being the position of the $n$-th vehicle, and $(x^{(n)},x^{(n-1)})$ taking the place of $(x,x^{(p)})$ in \eqref{E:OVM}) has been studied from a number of perspectives of stability (see \cite{chandler1958traffic,jiang2001full, kometani1958stability, wilson2011car}).  \emph{String stability} seeks to understand large-scale effects of stability in car-following dynamics by analysing a closed \lq\lq ring" of vehicles \cite{bando1995dynamical, gunter2020commercially}.  Starting with \eqref{E:OVM}, one can also build macroscopic models of traffic flow and congestion, and understand stability of such macroscopic behavior \cite{schonhof2007empirical,treiber2010three, tordeux2014collision}.  A important current research thrust is to understand various types of macrosopic behavior and stability when only a few vehicles in a platoon follow dynamics such as \eqref{E:OVM} \cite{cui2017stabilizing, stern2018dissipation,talebpour2015influence}; i.e., the effect of partial penetration of autonomous vehicles.

Our main focus here is the effect of \emph{noise} (and uncertainty) on car-following models.  The effect of random perturbations of traffic velocity upon road capacity has been (empirically) quantified for car-following models (viz. \cite{krauss1998microscopic}; cf. \cite{jost2003probabilistic}).  Ideas from physics of droplet clusters have also been used in modelling empirical traffic  observations \cite{mahnke2005probabilistic}.    Uhlenbeck-Ornstein processes have been used to model fluctuations (similar to the diffusive term in \eqref{E:OVM}) in driver behavior \cite{laval2014parsimonious}, explaining some instabilities and oscillations in observed data (see also \cite{yuan2018geometric}, and see \cite{ngoduy2019langevin} for an analysis of stability in the face of Cox-Ingersoll-Ross nonnegative perturbations).  See also \cite{matin2020nonlineara,matin2020nonlinearb}.  

\begin{figure}
    \centering
    \includegraphics[width =0.7\columnwidth]{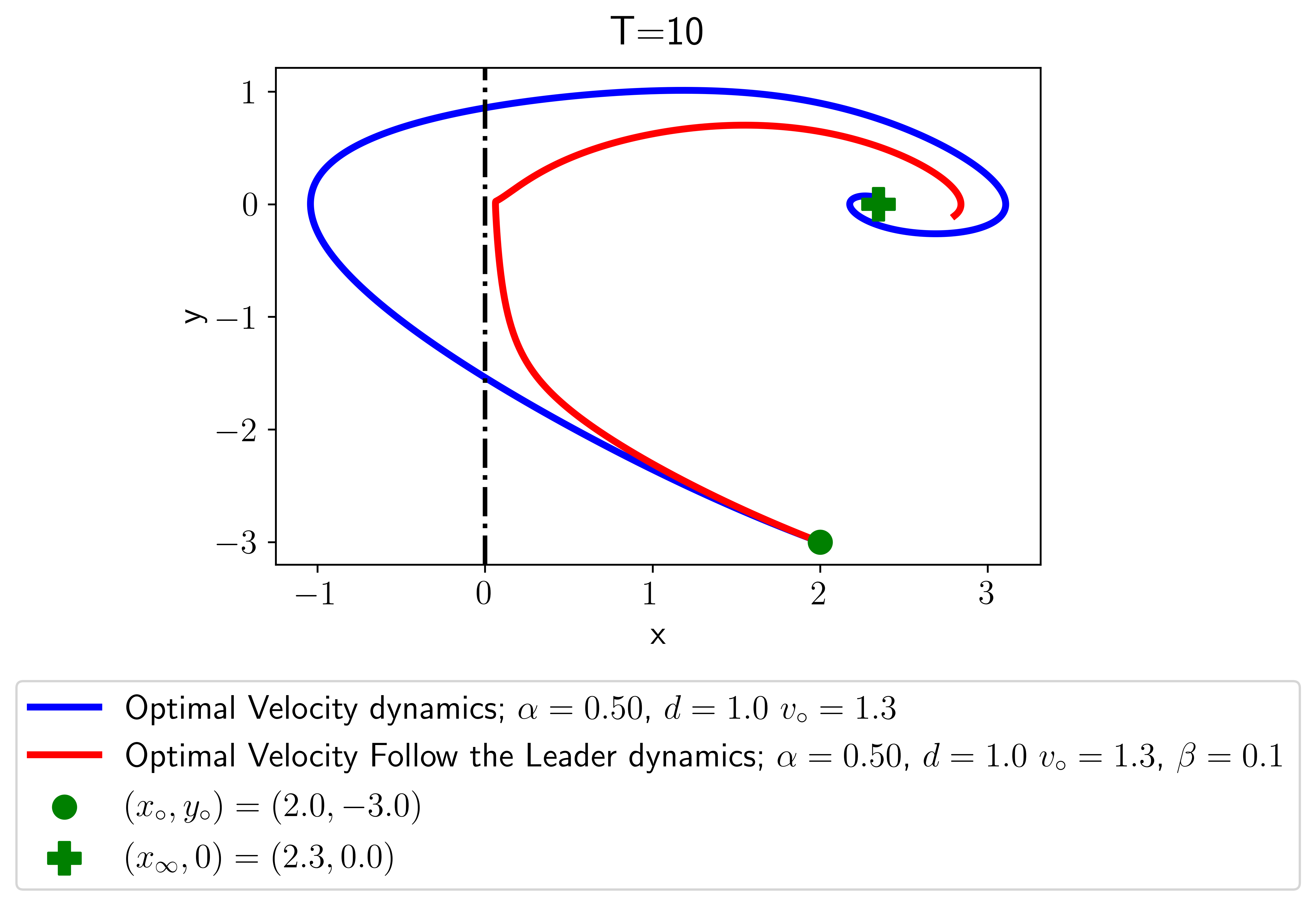}
    \caption{Solutions}
    \label{fig:solutions}
\end{figure}

In contrast to interests in stability (which focus on platoon behavior near equilibrium), our main interest in this paper is \emph{collisions}; where $x_t=x^{(p)}_t$ in \eqref{E:OVM}.  
If $\beta = \eps \equiv 0$ (the optimal velocity model of \cite{bando1994structure, bando1995dynamical}), the dynamics of \eqref{E:OVM} are continuous in a neighborhood of a collision and thus can be solved backward to identify states evolving to collision (see Figure \ref{fig:solutions}).

We will see (in Section \ref{S:deterministic_collision}) that the singularity (the $\beta$ term in \eqref{E:OVM}) in fact precludes collision in the absence of noise (we have not found a rigorous proof of this in the existing literature, so have included one here); as the following vehicle gets closer to the preceding vehicle, it decelerates quickly enough to prevent collision (see Figure \ref{fig:solutions}).

To simplify our focus even more, let's assume that the preceding vehicle has constant velocity $v_\circ$; let's assume that
\begin{equation*} x^{(p)}_t\Def v_\circ t \qquad t>0 \end{equation*}
for some $v_\circ$ in 
\begin{equation} \label{E:vrange} (0,V(\infty))=(0,1+\tanh(2))=(0,1.96) \end{equation}
i.e., the preceding starts at the origin and travels rightward with velocity $v_\circ$.   The restriction \eqref{E:vrange} allows us to invert $V$ at $v_\circ$ (see \eqref{E:equil}); $V(\infty)$ represents the maximum allowable velocity in our car-following model (abstractly, one can rescale time and space to ensure that $v_\circ$ is in the interval \eqref{E:vrange}, changing the constants $\alpha$, $\beta$, and $d$ in the process).
We also denote the initial position of the following vehicle by $-x_\circ$ with $x_\circ>0$; the following vehicle starts distance $x_\circ$ \lq\lq behind" the preceding vehicle.
Let's transform \eqref{E:OVM} by the Galilean change of coordinates
\begin{equation}\label{E:model_reduced} 
X^\eps_t  \Def v_\circ t - x_t \qquad \text{and}\qquad 
Y^\eps_t  \Def v_\circ - \dot x_t,\end{equation}
so $X^\eps$ and $Y^\eps$ are the position and velocity of the following vehicle behind the preceding vehicle (i.e., the following distance).
Then \eqref{E:OVM} becomes the $\R^2$-valued\footnote{We endow $\R^2$ with its standard linear operations and norm} Ito stochastic differential equation
\begin{equation}\label{E:initial_stochastic}\begin{aligned}
    d X^\eps_t &= Y^\eps_t dt\\
    dY^\eps_t  &=  \lb- \alpha \lb V \left(\frac{X^\eps_t}{d}\right)  - v_\circ +Y^\eps_t\rb - \beta \frac{Y^\eps_t}{(X^\eps_t)^2}\rb dt - \eps d W_t\\
    (X^\eps_0, Y^\eps_0) &= (x_\circ,y_\circ).
\end{aligned}\end{equation}
 Here 
\begin{equation} \label{E:ycircdef} y_\circ\Def  v_\circ-\dot x_0.\end{equation}
Tracing back to \eqref{E:model_reduced}, collision in \eqref{E:initial_stochastic} corresponds to $X^\eps_t=0$.
Fixing an underlying probability triple $(\Omega,\filt,\BP)$ (with associated expectation operator $\BE$), we assume that 
$W$ is a Brownian motion (on $(\Omega,\filt,\BP)$) with respect to a right-continuous filtration $\{\filt_t\}_{t\ge 0}$.

The main focus of this paper is \emph{near-collision dynamics}.  From an applied standpoint, a mathematical understanding of near-collision dynamics in \eqref{E:initial_stochastic} can give tools and a framework to understand various tradeoffs and optimal collision-avoidance strategies.  From a mathematical standpoint, we are interested in the interaction between noise and singularity in \eqref{E:initial_stochastic}; similar questions of interaction between noise and singularities are found in molecular dynamics \cite{morita1988brownian}, and typically involve careful analyses of boundary layers \cite{de1996theory,sowers2005boundary}.

Some prior works address questions of collision in deterministic models from different perspectives. For human-driven vehicles, \cite{davis2003modifications} uses simulation of the optimal velocity model (and several modifications) to study the effect of different values of driver's reaction time on likelihoods of collision.  Similarly, simulation of several microscopic acceleration models \cite{hamdar2008existing} have suggested quantifiable relationships between variations in safety conditions and macroscopic traffic patterns (in which collisions are allowed).  Finally, \cite{davis2014nonlinear} considers the dynamics of autonomous vehicles, and shows that, under the assumption of stability, collision does not happen for parameters in certain ranges.

Our efforts are organized as follows. In Section \ref{S:deterministic_collision}, we develop rigorous boundary-layer analysis for the deterministic system arising from \eqref{E:initial_stochastic} with $\eps=0$.  This depends on a transformation of coordinates near collision.  Our main result in this section, Theorem \ref{T:absorbing}, is that collisions cannot happen in \eqref{E:initial_stochastic} with $\eps=0$.  In Section \ref{S:PSPD}, we consider small-noise asymptotics of  \eqref{E:initial_stochastic} (i.e., $\eps \searrow 0$). A number of the techniques and results of Section \ref{S:PSPD} are modifications of the deterministic calculations of Section \ref{S:deterministic_collision}.  The main result of this section, Theorem \ref{T:main}, is that collision is unlikely over long time intervals for small $\eps$, with the length of the time interval growing as $\eps\searrow 0$.  Our arguments of Section \ref{S:PSPD} heavily depend on rigorous tools from stochastic analysis and our results represent provable bounds.  We furthermore hope that our arguments will generalize to a range of other types of singular terms in collision-avoidance models \cite{aghabayk2015state, gazis1961nonlinear}.

We note that a number of Lyapunov-based results have been proven away from the collision region (where one could compare behavior with a linearized system) \cite{tordeux2014collision, davis2014nonlinear}. Our results strongly depend on initial conditions of the following vehicle; our background interest is in the boundary layer near collision.  Our results can be modified to cover regimes where the following vehicle is asymptotically close to the lead vehicle (cf. \cite{davis2014nonlinear,tordeux2014collision}).

This work is based on the results of \cite{MatinThesis}.  

\section{Collisions in Deterministic Dynamics}\label{S:deterministic_collision}

Let's start by writing \eqref{E:initial_stochastic}, with $\eps=0$, in the framework of an ordinary differential equation.   Define the state space
\begin{equation}\label{E:statespace}  \Space\Def (0,\infty)\times \R. \end{equation}
The boundary 
\begin{equation} \label{E:collisionspace}\bdy \Space=\{0\}\times \R \end{equation}
where $x=0$ corresponds to collision.  Defining $\BB: \Space \to \R^2$ by 
\begin{equation}
    \label{E:B_dynamic}
    \BB(z) \Def \left(y,- \alpha \lb V\left(\frac {x}{d} \right) - v_\circ + y \rb  - \beta \frac{y}{x^2}\right), \qquad z=(x,y)\in \Space
\end{equation}
consider the ODE
\begin{equation}\label{E:OV_deterministic_delta_0}
    \dot z_t = \BB(z_t)
\end{equation}
starting at 
\begin{equation} \label{E:ZinitialDef} z_0= (x_0,y_0)=z_\circ\Def (x_\circ,y_\circ)\in \Space
\end{equation}
with $y_\circ$ as in \eqref{E:ycircdef}.
We write $z_t=(x_t,y_t)$; this is the dynamics of \eqref{E:initial_stochastic} with $\eps=0$.
Since $\BB$ is smooth on $\Space$, standard ODE theory ensures that \eqref{E:OV_deterministic_delta_0} has a unique solution on a maximal forward interval 
\begin{equation} \label{E:MIDDef} [0,\MID)\subset \R_+\end{equation}
of definition.

We are interested in understanding how collisions can occur in the model  \eqref{E:OV_deterministic_delta_0}.  In light of the state space \eqref{E:statespace} and its boundary \eqref{E:collisionspace}, questions of collision in \eqref{E:OV_deterministic_delta_0} correspond to questions of hitting the boundary $\bdy \Space$ of $\Space$ and the maximal time $\MID$. This motivates a careful boundary layer analysis of the dynamical model \eqref{E:OV_deterministic_delta_0}. 

We first of all understand the global structure of \eqref{E:OV_deterministic_delta_0}, writing it as a damped Hamiltonian system, thus giving us some macroscopic bounds.  We then show that $\MID$ in \eqref{E:MIDDef} must in fact be infinite (precluding collision in finite time).  We do this by assuming the converse (that $\MID<\infty$) and then precluding a number of ways in which collision could occur.  In the final case, we have to carefully understand the effect of the singularity in \eqref{E:OV_deterministic_delta_0} at $x=0$.  We analyze this case by understanding the invariant manifold of a solution of \eqref{E:OV_deterministic_delta_0} and showing that it cannot intersect $\bdy \Space$.  We then additionally show global asymptotic stability by an argument involving the Poincar\'e-Bendixon theorem,

Let's start by understanding some global properties of \eqref{E:OV_deterministic_delta_0}.
Since $V:(0,\infty)\to (0,V(\infty))$ is an invertible transformation, the dynamics \eqref{E:OV_deterministic_delta_0} has a unique equilibrium point at  $(x_\infty, 0)\in \Space$, where
\begin{equation}\label{E:equil}
    x_\infty = d \cdot V^{-1}(v_\circ)=d\lb 2 + \tanh^{-1}(v_\circ + \tanh(-2))\rb,
\end{equation}
which is well-defined under the assumption that $v_\circ$ is in \eqref{E:vrange}. Namely, if the preceding vehicle has constant velocity $v_\circ$, the system \eqref{E:initial_stochastic} (with $\eps=0$, or alternately \eqref{E:OV_deterministic_delta_0}) selects $x_\infty$ as the proper following distance (and does so with velocity $v_\circ$).

\begin{figure}
    \centering
    \includegraphics[width=0.45\columnwidth]{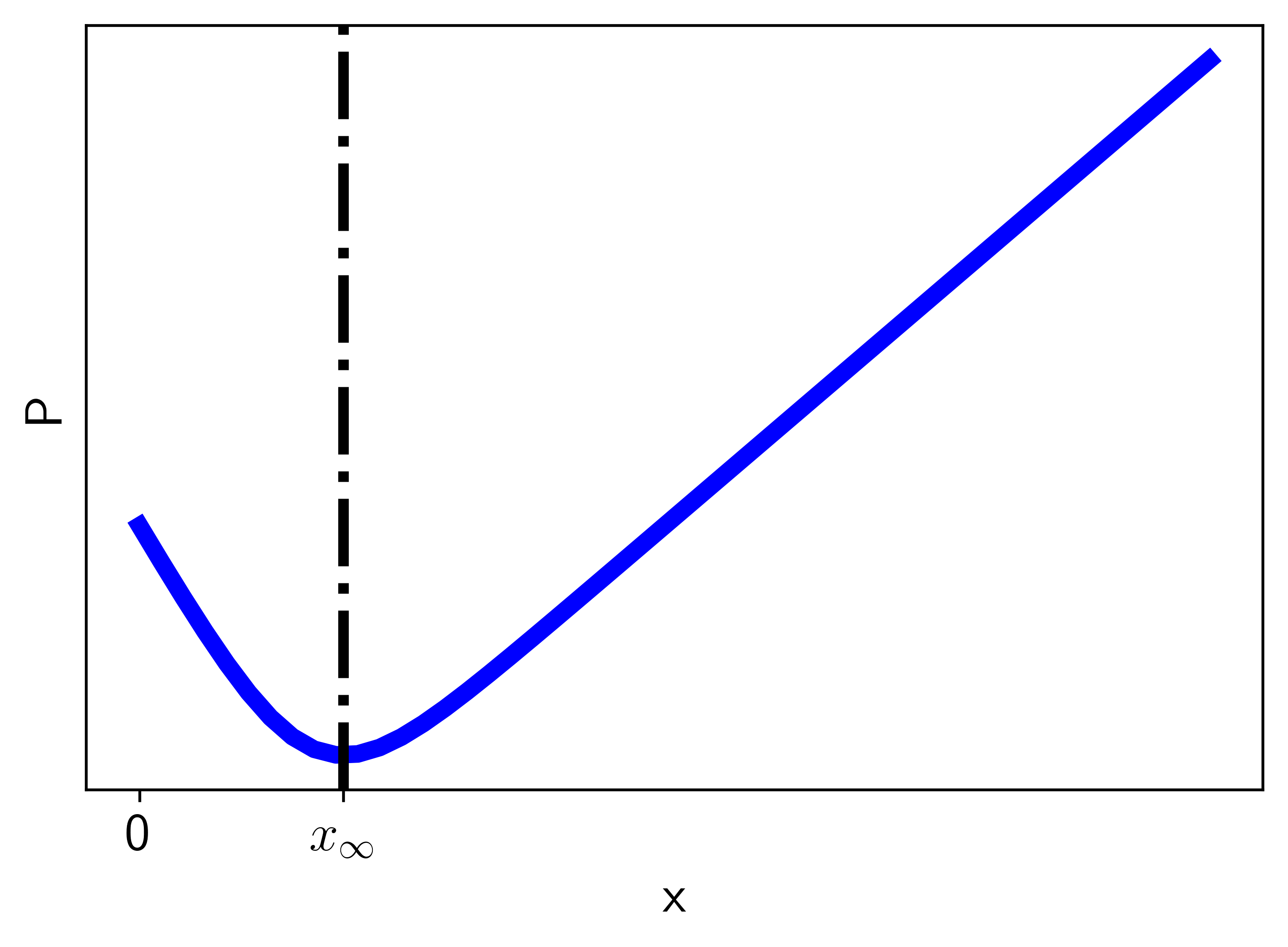}
    \includegraphics[width=0.45\columnwidth]{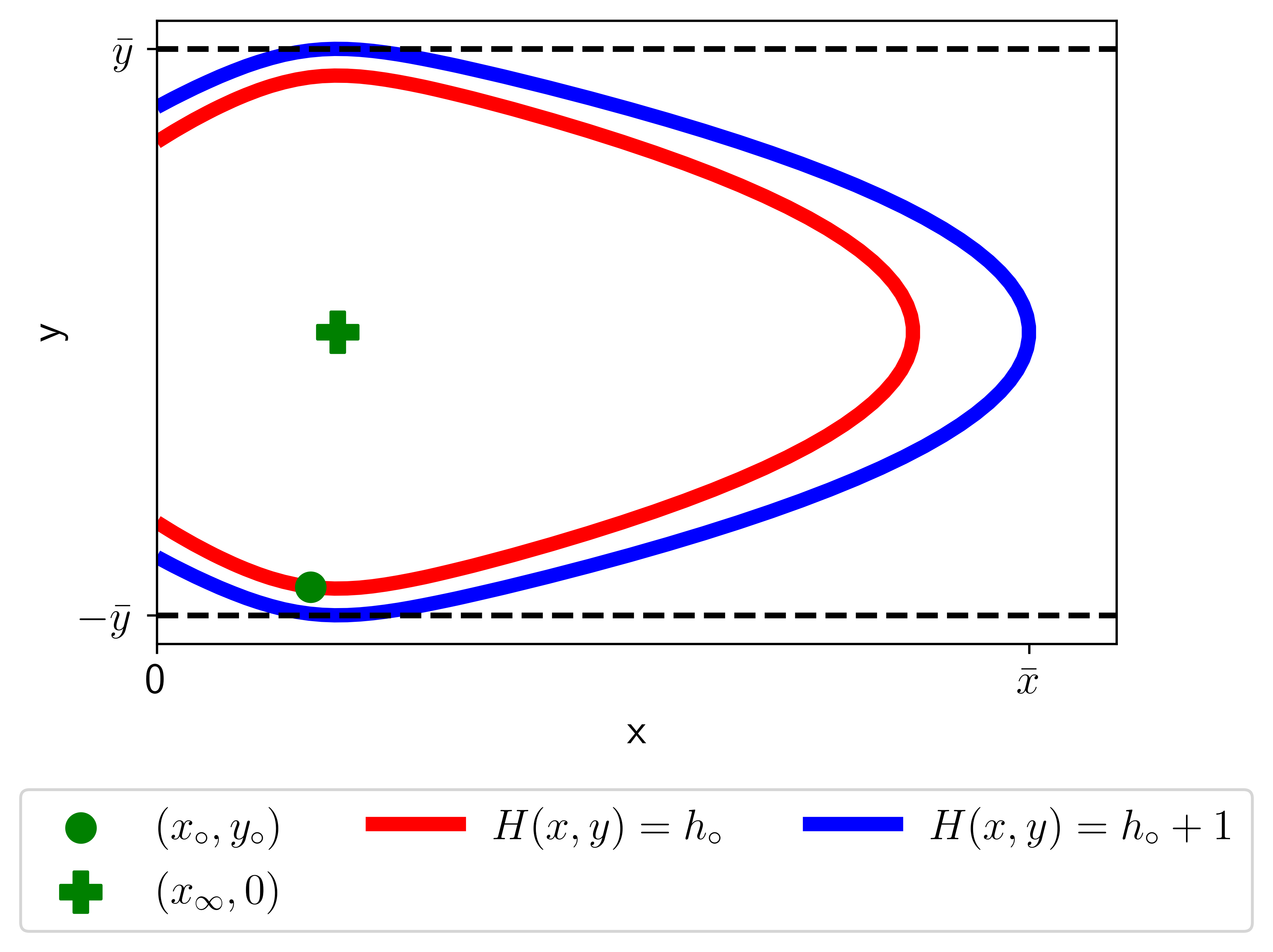}
    \caption{Potential $P$ and contour plot of $H$}
    \label{fig:geometry}
\end{figure}

Let's next consider a more global understanding of the dynamics of \eqref{E:OV_deterministic_delta_0}.  Define
\begin{equation*}
    H(x, y) \Def \frac 12 y^2 + P(x) \qquad (x, y) \in \R^2
\end{equation*}
where
\begin{equation*}
    P(x) \Def \alpha \int_{x'=x_\infty }^x \lb V \left(\frac{x'}{d} \right) - v_\circ \rb dx'; \qquad x \in \R
\end{equation*}
see Figure \ref{fig:geometry}.  Since 
\begin{equation*}
    P'(x) = \alpha \lb V \left( \frac xd \right) - v_\circ \rb, \qquad x>0
\end{equation*}
we can write the dynamics \eqref{E:OV_deterministic_delta_0} as a \emph{damped Hamiltonian system}
\begin{equation}\label{E:main_deterministic} \begin{aligned} \dot x_t&=y_t=\frac{\partial H}{\partial y}(x_t,y_t) \\
\dot y_t &= -P'(x_t)-\alpha y_t -\beta\frac{y_t}{x_t^2} = -\frac{\partial H}{\partial x}(x_t,y_t)-\lb \alpha +  \frac{\beta}{x_t^2}\rb y_t \end{aligned} \end{equation}
for $t\in [0, \MID)$.

\noindent We can explicitly calculate
\begin{equation*}
        P(x) 
        = - \frac{\alpha d}{2} \log \left( \frac{1 - \tanh^2\left(2 - \frac{x}{d}\right)}{1 - \tanh^2\left(2 - \frac{x_\infty}{d}\right)} \right)  + \alpha \left(\tanh(2) - v_\circ \right) (x - x_\infty).
\end{equation*}
\begin{remark}\label{R:Pstructure} The monotonicity of $V$ implies that $V(x/d)-v_\circ <0$ if $x < x_\infty$ and $V(x/d)-v_\circ > 0$ if $x>x_\infty$. Thus $P$ is increasing on $(x_\infty,\infty)$, decreasing on interval $(0,x_\infty)$ and has a minimum at $x_\infty$.  Furthermore $\lim_{x\nearrow \infty}P(x)=\infty$. \end{remark}
\noindent We then have
\begin{proposition}\label{P:decreasing}
The point $(x_\infty,0)$ is the unique equilibrium point of \eqref{E:OV_deterministic_delta_0}.  Secondly, the function $H$ is nonincreasing along the trajectory of \eqref{E:OV_deterministic_delta_0}, and strictly decreasing if and only if $(x_\circ,y_\circ)\not =(x_\infty,0)$
\end{proposition}
\begin{proof} 
An explicit calculation shows that $(x_\infty,0)$ is indeed the unique equilibrium point of \eqref{E:OV_deterministic_delta_0}.

For $t\in [0,\MID)$,
\begin{equation}\label{E:energyLevel}
    \frac{dH}{dt}(x_t, y_t)= \dot x_t \frac{\partial H}{\partial x}(x_t, y_t) + \dot y_t \frac{\partial H}{\partial y}(x_t, y_t) = -\lb \alpha + \frac{\beta}{x_t^2}\rb y_t^2 \le 0
\end{equation}
so $H$ is nonincreasing.

Assume now that $H$ is strictly decreasing along trajectories of \eqref{E:main_deterministic}.  This is impossible if $(x_\circ,y_\circ)$ is the fixed point $(x_\infty,0)$, so necessarily $(x_\circ,y_\circ)\not = (x_\infty,0)$.  On the other hand, assume that $(x_\circ,y_\circ)\not = (x_\infty,0)$, and then assume that $H$ is not strictly decreasing.  Then there is a nonempty $(t_1,t_2)\subset (0,\infty)$ on which $H(x_t,y_t)$ is constant for $t\in (t_1,t_2)$.  From \eqref{E:energyLevel}, this means that $y_t=0$ for $t\in (t_1,t_2)$.  From \eqref{E:main_deterministic}, we then have that $x$ takes on a constant value $x_c$ on $(t_1,t_2)$.  This means that $(x_c,0)$ is a fixed point and thus must be the unique fixed point $(x_\infty,0)$.  This violates our assumption, so necessarily $H$ must be strictly decreasing along trajectories of \eqref{E:main_deterministic}.
\end{proof}
\noindent Geometrically, Proposition \ref{P:decreasing} implies that the sublevel sets of $H$ are invariant under the dynamics of  \eqref{E:OV_deterministic_delta_0}.   Define 
\begin{equation}\label{E:hcircdef} h_\circ \Def H(x_\circ,y_\circ) \end{equation}
which is the initial value of $H(x_t,y_t)$.  Define
\begin{equation*}   \bar y \Def \sqrt{2 (h_\circ+1)}\qquad \text{and}\qquad 
\bar x \Def \left(P\Big|_{(x_\infty, \infty)}\right)^{-1}(h_\circ+1) \end{equation*}
and fix $(x,y)\in \Space$.  Recall Remark \ref{R:Pstructure}.  If $x> \bar x$, then the fact that $P$ is increasing on $(x_\infty,\infty)$ implies that
\begin{equation*} H(x,y)\ge P(x)>P(\bar x)=h_\circ+1 \end{equation*}
Similarly, if $|y|>\bar y$, the nonnegativity of $P$ allows us to write
\begin{equation*} H(x,y)\ge \tfrac12 y^2> \tfrac12 \bar y^2= h_\circ+1. \end{equation*}
Thus
\begin{equation}\label{E:levelsetscompact} \lb (x,y)\in \Space: H(x,y)\le h_\circ\rb \subset  \lb (x,y)\in \Space: H(x,y)\le h_\circ+1\rb \subset (0,\bar x]\times [-\bar y,\bar y] \end{equation}
(we will here only use the fact that the leftmost set of \eqref{E:levelsetscompact} is contained in the rightmost set; in Section \ref{S:PSPD}, the extra margin of $h_\circ+1$, as opposed to $h_\circ$, will give us a margin to allow for noise).  See Figure \ref{fig:geometry}.

Proposition \ref{P:decreasing} hence implies that
\begin{equation}\label{E:compactset} (x_t,y_t)\in (0,\bar x]\times [-\bar y,\bar y] \end{equation}
for all $0\le t<\MID$.  

We can now organize some standard calculations to understand a bit more about collision.
\begin{lemma}\label{L:collision} If $\MID<\infty$, \begin{equation} \label{E:collisionexit}\lim_{t\nearrow \MID}x_t=0.\end{equation}
\end{lemma}
\begin{proof} 
Fix an increasing sequence $(t_n)_{n\in \N}$ of times in $[0,\MID)$ such that $\lim_{n\nearrow \infty}t_n=\MID$.  Keeping \eqref{E:compactset} in mind and using the fact that $[0,\bar x]\times [-\bar y,\bar y]$ is a compact subset of $\R^2$, we can find a subsequence $(t_{n_k})_{k\in \N}$ of $(t_n)_{n\in \N}$ such that $(x^*,y^*)=\lim_{k\to \infty} (x_{t_{n_k}},y_{t_{n_k}})$ exists (as a limit in $\R^2$); it also follows that $(x^*,y^*)$ is in the closure $\bar \Space$ of $\Space$.  If $(x^*,y^*)\in \Space$, a solution of \eqref{E:OV_deterministic_delta_0} can be constructed from $(x^*,y^*)$ backward and forward in time.  By the uniqueness of the backward solution, we in fact have that $\lim_{t\nearrow \MID}(x_t,y_t)=(x^*,y^*)$.  The existence of the forward solution contradicts the definition of $\MID$; this implies that $(x^*,y^*)\not \in \Space$, so in fact $(x^*,y^*)\in \bar \Space\setminus \Space = \bdy \Space$.

To summarize, every increasing sequence $(t_n)_{n\in N}$ converging to $\MID$ has a subsequence $(t_{n_k})_{k\in \N}$ such that $\lim_{k\nearrow \infty}x_{t_{n_k}}=0$.  The claim follows.
\end{proof}
\noindent The asymptotics of \eqref{E:collisionexit} amount to collision in finite time; recall \eqref{E:collisionspace}.
A boundary-region analysis near $\bdy \Space$ will preclude this case.  The main result of this section is
\begin{theorem}\label{T:absorbing} We have that $\MID=\infty$ and the equilibrium solution $(x_t,y_t) \equiv (x_\infty, 0)$ \textup{(}defined in \eqref{E:equil}\textup{)} of the dynamical model \eqref{E:OV_deterministic_delta_0} is globally asymptotically stable.
\end{theorem}
\noindent We will formalize the proof of this at the end of the section.

For now (prior to the proof of Theorem \ref{T:absorbing}), let's assume that
\begin{equation} \label{E:Tfinite} \boxed{\boxed{\MID<\infty;}} \end{equation}
We will show that this leads to a contradiction.

Let's also fix
\begin{equation}\label{E:xminusdef} x_-\Def \min\{x_\circ,x_\infty\}, \end{equation}
(with $x_\infty$ as in \eqref{E:equil}); we are interested in the dynamics of \eqref{E:OV_deterministic_delta_0} in the boundary region
\begin{equation*}(0,x_-)\times \R \end{equation*}
to the left of the initial condition and the fixed point $(x_\infty,0)$.
Define 
\begin{equation}\label{E:checkt} \check t\Def \sup\lb t<\MID: x_t\ge \tfrac12 x_-\rb; \end{equation}
since $x_\circ>\tfrac12 x_-$ (by \eqref{E:xminusdef}), assumption \eqref{E:Tfinite} and Lemma \ref{L:collision} imply that $\check t$ is well-defined.  Informally, $\check x_t$ is the last time that the solution of \eqref{E:OV_deterministic_delta_0} is to the right of $\tfrac12 x_-$.  

\begin{figure}
    \centering
    \includegraphics[width=0.5\columnwidth]{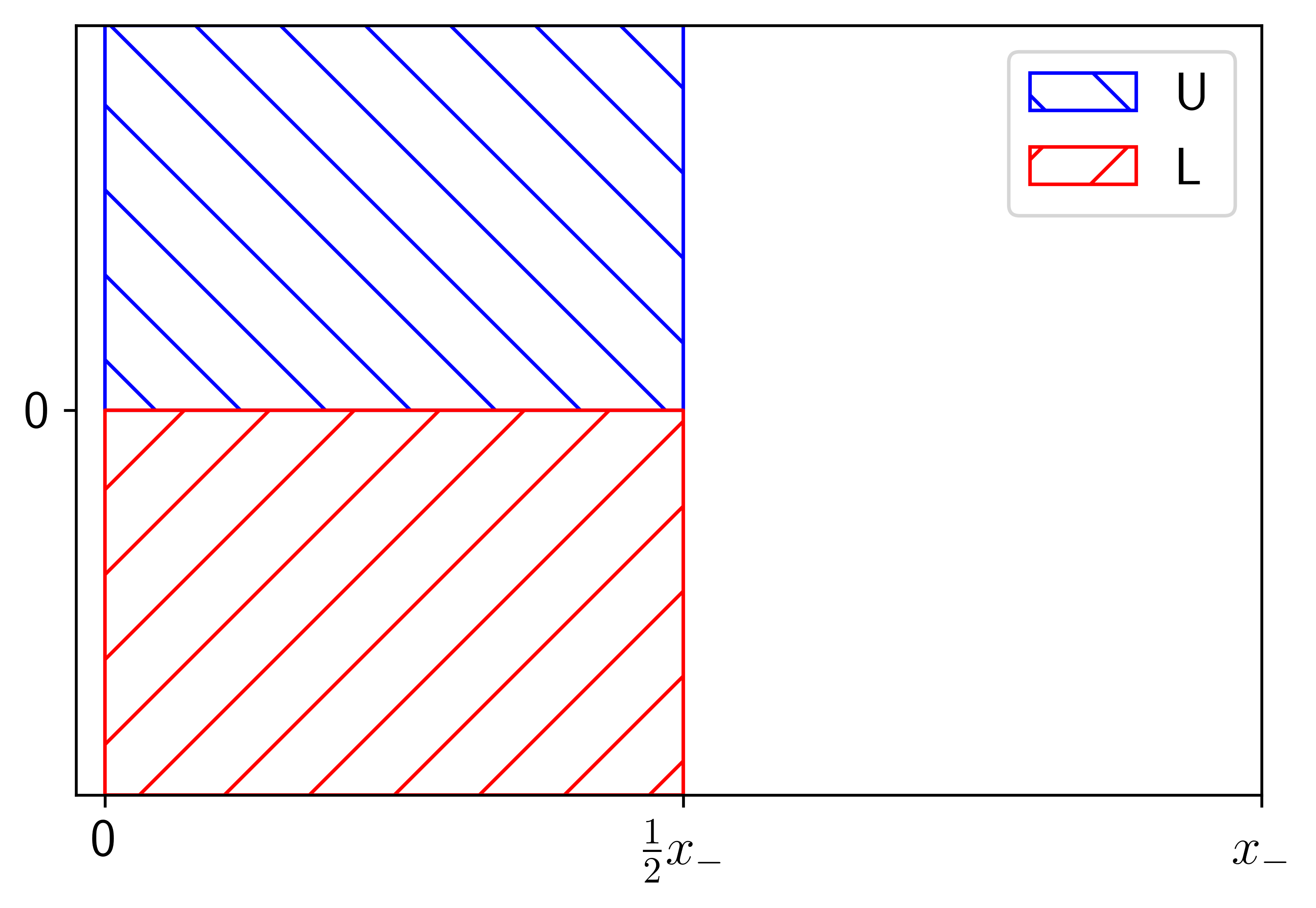}
    \caption{Boundary regions for deterministic dynamics}
    \label{fig:RUL}
\end{figure}
Let's write
\begin{equation*}\Region = \Lower \cup \Upper \end{equation*}
where
\begin{equation*} \Region\Def \left(0,\tfrac12 x_-\right)\times \R, \qquad
\Lower\Def \left(0,\tfrac 12 x_-\right)\times (-\infty,0), \qquad \text{and}\qquad
\Upper\Def \left(0,\tfrac 12 x_-\right)\times [0,\infty), \end{equation*}
partitioning $\left(0,\tfrac 12 x_-\right)\times \R$ into \lq\lq upper" and \lq\lq lower" parts (see Figure \ref{fig:RUL}).  By definition \eqref{E:checkt} of $\check t$, 
\begin{equation}\label{E:region} (x_t,y_t)\in \Region \qquad t\in (\check t,\MID)\end{equation}
We have
\begin{lemma}\label{L:upperinvariant} The upper set $\Upper$ is invariant for $\{(x_t,y_t);\, t\in (\check t,\MID)\}$; if $(x_{t'},y_{t'})\in \Upper$ for some $t'\in (\check t,\MID)$, then $(x_t,y_t)\in \Upper$ for $t\in [t',\MID)$. \end{lemma}
\begin{proof}
Assume not; then there is a nonempty $\left[t^{\prime },t^{\prime\prime}\right)\subset (\check t,\MID)$ such that $(x_{t^\prime},y_{t^\prime})\in \Upper$ but $(x_s,y_s)\not \in \Upper$ for $s\in \left(t^\prime,t^{\prime\prime}\right)$.  By \eqref{E:region}, we then have that $(x_s,y_s)\in \Region\setminus \Upper=\Lower$ for $s\in \left(t^\prime,t^{\prime\prime}\right)$.
Letting $s\searrow t^\prime$, we thus have that
\begin{equation*} (x_{t^\prime},y_{t^\prime}) \in \Upper \cap \bar \Lower = \left(0,\tfrac 12 x_-\right)\times \{0\}. \end{equation*}
By definition of $\Lower$, $y_t<0$ for $t\in \left(t^\prime,t^{\prime\prime}\right)$, so by continuity $\dot y_{t^\prime}\le 0$.
However, \eqref{E:OV_deterministic_delta_0}
 implies that
 \begin{equation*} \dot y_{t^\prime}=-\alpha\lb V\left(\frac{x_{t^\prime}}{d}\right)-v_\circ\rb >0 \end{equation*}
 since $x_{t^\prime}<x_\infty$.  This is a contradiction, and the claim thus holds.
\end{proof}
\noindent We also have
\begin{lemma}\label{L:upperdriftright} If $(x_{t'},y_{t'})\in \Upper$, then $x_t\ge x_{t'}$ for $t\in (t',\MID)$. \end{lemma}
\begin{proof}  
Assuming that $(x_{t'},y_{t'})\in \Upper$, Lemma \ref{L:upperinvariant} implies that $(x_t,y_t)\in \Upper$ for $t\in (t',\MID)$, so
\begin{equation*} x_t=x_{t'}+\int_{s=t'}^{t}y_sds\ge x_{t'} \end{equation*}
for $t\in (t',\MID)$.
\end{proof}

Taken together, we get
\begin{lemma} \label{L:bottombox} We have that
\begin{equation}\label{E:bottombox} (x_t,y_t)\in \Lower \end{equation}
for all $t\in [\check t,\MID)$, and 
\begin{equation}\label{E:ynegative} y_{\MID_-}\Def \varlimsup_{t\nearrow \MID}y_t\le 0. \end{equation}
\end{lemma}
\begin{proof} If $(x_t,y_t)\in \Upper$, Lemmas \ref{L:upperinvariant} and \ref{L:upperdriftright} contradict the fact that $\lim_{t\nearrow \MID}x_t=0$ (Lemma \ref{L:collision}).  Equation \eqref{E:ynegative} follows from \eqref{E:bottombox}.
\end{proof}

Furthermore,
\begin{lemma}\label{L:yinc} The map $t\mapsto y_t$ has strictly positive derivative on $(\check t,\MID)$.\end{lemma}
\begin{proof}
If $(x_t,y_t)\in \Lower$, then $P'(x_t)<0$ and by \eqref{E:main_deterministic} we have
\begin{equation*} \dot y_t=-P'(x_t)-\lb \alpha+ \frac{\beta}{x_t^2}\rb y_t> 0 \end{equation*}
and the result follows.
\end{proof}

Lemma \ref{L:yinc} implies that 
\begin{equation*} y_{\check t}=y_{\MID-}-\int_{s=\check t}^{\MID}\dot y_sds<y_{\MID-}\le 0. \end{equation*}
Collecting things together, we must in fact have that
\begin{equation*} \lim_{t\nearrow \MID}(x_t,y_t)=(0,y^*) \end{equation*}
for some $y_*\in (y_{\check t},0]$, with $y$ strictly increasing on $(\check t,\MID)$.  

Let's now try to parametrize the integral curve $t \in (\check t, \mathcal T) \mapsto (x_t, y_t)$.  In view of Lemma \ref{L:yinc}, we should be able to locally write
\begin{equation}\label{E:local_parametrization}
    x_t = \varphi(y_t), \quad t \in (\check t, \MID)
\end{equation}
for some function $\varphi$ (which we will in fact construct in a moment),
and if so and $\varphi$ is differentiable,
\begin{equation}\label{E:localpara} \begin{aligned}
    y_t=\dot x_t &= \varphi'(y_t) \dot y_t= \varphi'(y_t) \lb - \alpha \lb V \left(\frac{x_t}{d}\right) - v_\circ + y_t \rb - \beta \frac{y_t}{x_t^2} \rb\\
        &= \varphi'(y_t)  \lb - \alpha \lb V \left(\frac{\varphi(y_t)}{d}\right)- v_\circ + y_t\rb - \beta \frac{y_t}{\varphi^2(y_t)} \rb.\end{aligned}
\end{equation}
Let's formalize this by defining
\begin{equation*}
    \ff(y, \Phi) \Def \frac{-y}{\alpha \lb V\left(\frac{\Phi}{d}\right)- v_\circ + y\rb + \beta \frac{y}{\Phi^2}}=\frac{-y}{P'(\Phi)+\{\alpha+\beta/\Phi^2\}y}. \qquad (y,\Phi)\in (-\infty, 0) \times (0, x_-)
\end{equation*}
Since $P'<0$ on $(0,x_-)$ and $\alpha$ and $\beta$ are also positive, $\ff$ is well-defined.

Let's now again appeal to the abstract theory of solutions of ODE's.
Since $(y_{\check t},x_{\check t})\in (-\infty,0)\times (0,x_-)$, there is a maximal interval $(\yy_-,\yy_+)\subset(-\infty,0)$ containing $y_{\check t}$ such that the ODE
\begin{equation}\label{E:phi_dynamic_deterministic}\begin{aligned} \varphi'(y)&=\ff(y,\varphi(y)) \\
\varphi(y_{\check t})&=x_{\check t} \end{aligned}\end{equation}
has a unique solution on $(\yy_-,\yy_+)$.  By \eqref{E:checkt}, we have that $x_{\check t}=\tfrac12 x_-$.  By definition of the domain of $\ff$, $\varphi(y) \in (0, x_-)$ for $y\in (\yy_-,\yy_+)$ and thus $P'(\varphi(y))< 0$ for $y\in (\yy_-, \yy_+)$.  Since $(\yy_-,\yy_+)\subset(-\infty,0)$, $\ff(y,\varphi(y))<0$ for $y\in (\yy_-,\yy_+)$; thus $\varphi$ is strictly decreasing on $(\yy_-,\yy_+)$.

\begin{lemma}\label{L:collision_avoid}  We have that
\begin{equation*} \inf_{y\in [y_{\check t},\yy_+)}\varphi(y)>0. \end{equation*}
\end{lemma}
\begin{proof}
Let's define a reference curve
\begin{equation*} R(y)\Def \lb \left(\tfrac12 x_-\right)^{-1}+\frac{1}{\beta}(y-y_{\check t})\rb^{-1}
\end{equation*}
for $y\in [y_{\check t},0)$.  $R(y_{\check t})=\tfrac12 x_-=\varphi(y_{\check t})$; we will see that $R$ acts as a lower bound for $\varphi$ on $[y_{\check t},\yy_+)$.  We also have that
\begin{equation*} \frac{\dot R(y)}{R^2(y)}=-\frac{1}{\beta} \end{equation*}
for $y\in (y_{\check t},0)$.

Let's now define
\begin{equation*} \xi(y)\Def \left(\frac{1}{\varphi(y)}-\frac{1}{R(y)}\right)^+=\left( \frac{1}{\varphi(y)}-\frac{1}{R(y)}\right)\bOne_{\{\nicefrac{1}{\varphi(y)}>\nicefrac{1}{R(y)}\}}=\left( \frac{1}{\varphi(y)}-\frac{1}{R(y)}\right)\bOne_{\{\varphi(y)<R(y)\}} \end{equation*}
for $y\in [y_{\check t},\yy_+)$, which measures the amount by which $1/\varphi$ is larger than $1/R$.
If $y\in (y_{\check t},\yy_+)$ and $\varphi(y)< R(y)$, then
\begin{align*} \dot \xi(y) &= \left( \frac{\dot R(y)}{R^2(y)}-\frac{ \varphi'(y)}{\varphi^2(y)}\right)
=\left(-\frac{1}{\beta}-\frac{-y}{\alpha \varphi^2(y) \lb V\left(\frac{\varphi(y)}{d}\right)- v_\circ+ y \rb + \beta y}\right)\\
&=\left(\frac{(-y)}{\beta (-y)+\alpha \varphi^2(y)(-y)-\alpha \varphi^2(y) \lb V\left(\frac{\varphi^2(y)}{d}\right)- v_\circ\rb}-\frac{1}{\beta}\right).
\end{align*}
Since $\varphi$ is decreasing on $[y_{\check t},\yy_+)$, 
\begin{equation*} \varphi(y)\le \varphi(y_{\check t})=x_{\check t}=\tfrac12 x_-<x_\infty \end{equation*}
for $y\in [y_{\check t},\yy_+)$, so 
\begin{equation*} V\left(\frac{\varphi(y)}{d}\right)-v_\circ\le V\left(\frac{x_\infty}{d}\right)-v_\circ=0 \end{equation*}
for $y\in [y_{\check t},\yy_+)$.  Again using the fact that $[y_{\check t},\yy_+)\subset (-\infty,0)$, we have that
\begin{equation*} \alpha \varphi^2(y)(-y)-\alpha \varphi^2(y) \lb V\left(\frac{\varphi(y)}{d}\right)- v_\circ\rb>0 \end{equation*} and hence
\begin{equation*} \frac{(-y)}{\beta (-y)+\alpha \varphi^2(y)(-y)-\alpha \varphi^2(y) \lb V\left(\frac{\varphi^2(y)}{d}\right)- v_\circ\rb}<\frac{1}{\beta}\end{equation*}
for $y\in [y_{\check t},\yy_+)$.
Hence
\begin{equation*} \dot \xi(y)\le 0 \end{equation*}
if $y\in [y_{\check t},\yy_+)$ and $\varphi(y)<R(y)$.  
Thus $\xi$ is decreasing on $[y_{\check t},\yy_+)$.  Since
\begin{equation*} \xi(y_{\check t})=\left(\frac{1}{x_{\check t}}-\frac{1}{x_{\check t}}\right)^+=0, \end{equation*}
we in fact have that $\xi\le 0$ on $[y_{\check t},\yy_+)$.
Thus $\varphi\ge R$ on $[y_{\check t},\yy_+)$ and consequently
\begin{equation*} \inf_{[y_{\check t},\yy_+)}\varphi(y)
\ge \inf_{[y_{\check t},0)}R(y)\ge \lb \left(\tfrac12 x_-\right)^{-1}+\frac{1}{\beta}(-y_{\check t})\rb^{-1}>0, \end{equation*}
giving the claim.
 \end{proof}

We can now verify that collision is impossible in the deterministic model of \eqref{E:OV_deterministic_delta_0}.
\begin{proof}[Proof of Theorem \ref{T:absorbing}]
Uniqueness allows us to formalize \eqref{E:local_parametrization}.  
In light of Lemma \ref{L:yinc}, $y$ is in fact a bijection from $[\check t,\MID)$ to $[y_{\check t},y_{\MID-})$ with continuous inverse.  Inverting this on the interval $[y_{\check t},\yy_+)\cap [y_{\check t},y_{\MID-})$ and noting that the continuous image of a connected set is connected, we get that
\begin{equation*} y\big|_{[\check t,\MID)}^{-1}\left([y_{\check t},\yy_+)\cap [y_{\check t},y_{\MID-})\right)=[\check t,\MID') \end{equation*}
for some $\MID'\le \MID$.  From \eqref{E:localpara}, 
\begin{equation*} (\varphi(y_t),y_t) \end{equation*}
satisfies \eqref{E:OV_deterministic_delta_0} for $t\in [\check t,\MID')$ and thus must agree with $(x_t,y_t)$ in $[\check t,\MID')$.

We will show that our standing assumption of \eqref{E:Tfinite} leads to a contradiction.  Assume first that $\MID'<\MID$.  This then implies that $(x_{\MID'},y_{\MID'})$ exists and is in fact equal to $(x_{\MID'-},y_{\MID'-})$.  By definition of $\MID'$, we must then have that $(x_{\MID'-},y_{\MID'-})\in \bdy \left(\left(0, x_-\right)\times (-\infty,0)\right)$.  Since $\MID'<\MID$, we must also have that $(x_{\MID'},y_{\MID'})\in \Space$.  By Lemma \ref{L:bottombox}, we must also have that
$(x_{\MID'},y_{\MID'})\in \Lower$.  However,
\begin{equation*} \bdy \left(\left(0, x_-\right)\times (-\infty,0)\right) \cap \Lower=\emptyset \end{equation*}
which is a contradiction.  On the other hand, if $\MID'=\MID$, then 
\begin{equation*} \lim_{t\nearrow \MID'}\varphi(y_t)=\lim_{t\nearrow \MID}x_t=0, \end{equation*}
which is precluded by Lemma \ref{L:collision_avoid}.

To show the global and asymptotic stability of the equilibrium solution, recall \eqref{E:compactset}.  The orbit of the flow $(x_t, y_t)$ is forward complete and contained in $[0,\bar x]\times [-\bar y,\bar y]$.  In light of the Poincar\'e-Bendixson theorem \cite{MR0069338}, the $\omega$-limit set of $(x_t,y_t)$ is either a fixed point or a periodic orbit.  By Proposition \ref{P:decreasing}, $(x_\infty,0)$ is the unique fixed point.  The function $H$ must be periodic along periodic orbits, but if $(x_\circ,y_\circ)\not = (x_\infty,0)$, $t\mapsto H(x_t,y_t)$ must be strictly decreasing on $[0,\MID)$ (again, by Proposition \ref{P:decreasing}).  This precludes the existence of periodic orbits, implying the claimed asymptotic stability.
\end{proof}

\section{Collisions in Stochastic Dynamics}\label{S:PSPD}

\begin{figure}
    \centering
    \includegraphics[width =0.7\columnwidth]{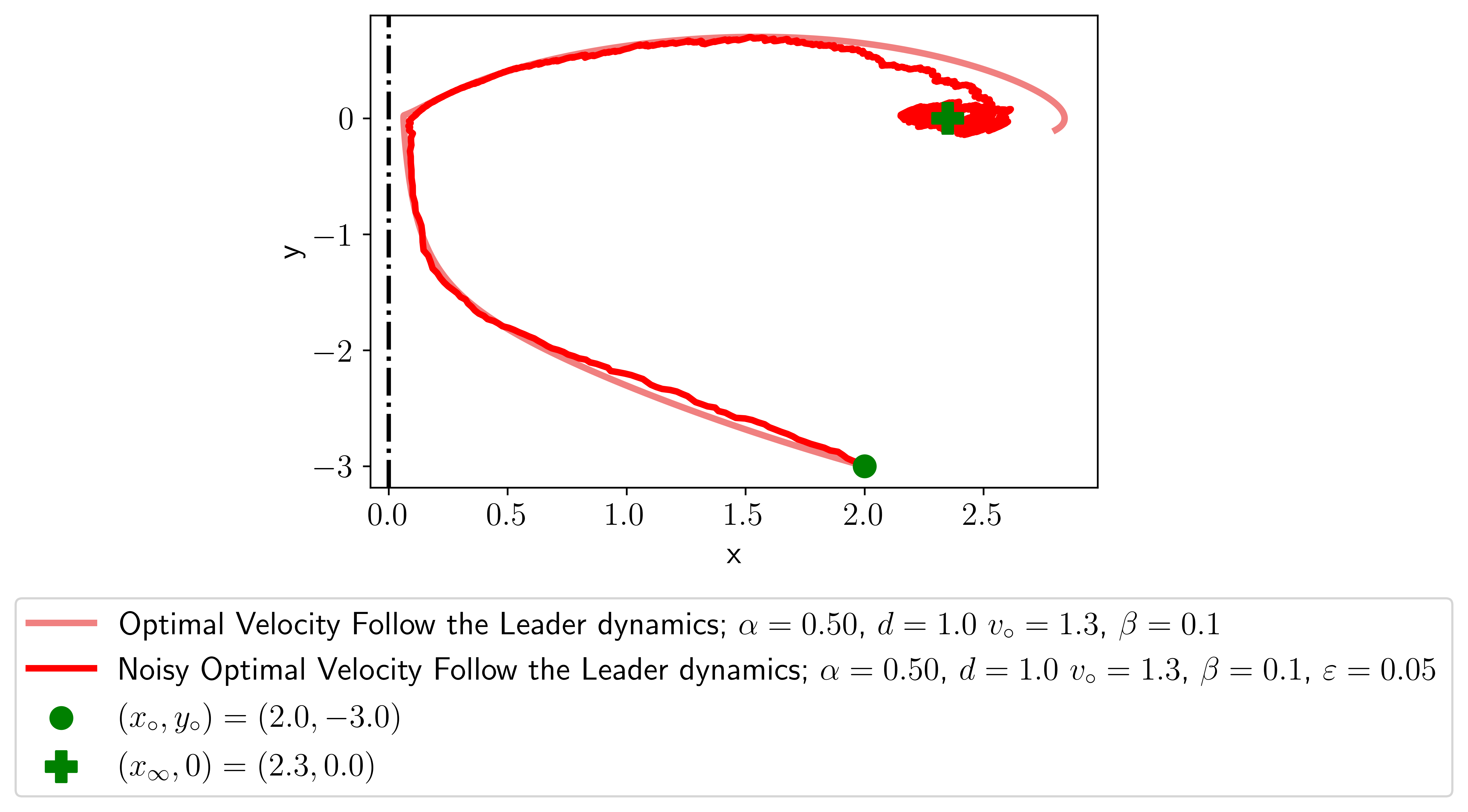}
    \caption{Noisy perturbations of dynamics}
    \label{fig:solutionswithnoise}
\end{figure}

We next consider the effect of small stochasticity in \eqref{E:initial_stochastic}. In the deterministic case, we showed that the system \eqref{E:OV_deterministic_delta_0} never hits the collision boundary \eqref{E:collisionspace} (more precisely, we showed that the maximal interval of definition of \eqref{E:MIDDef} was all of $\R_+$).
We want to prove similar results here. We will first of all understand how to construct the solution of \eqref{E:initial_stochastic}; we will carefully regularize the singularity in \eqref{E:initial_stochastic} and then use standard tools from the theory of stochastic differential equations to construct approximate solutions.
This will allow us to formalize a usable definition of the solution of \eqref{E:initial_stochastic}, which allow us to prove global bounds on the solution of \eqref{E:initial_stochastic}.  We will then prove the analogue of Theorem \ref{T:absorbing}, showing that collision is unlikely over long time intervals (where the length of the time interval grows as size of the noise becomes small).  See Figure \ref{fig:solutionswithnoise}

\begin{figure}
    \centering
    \includegraphics[width=0.5\columnwidth]{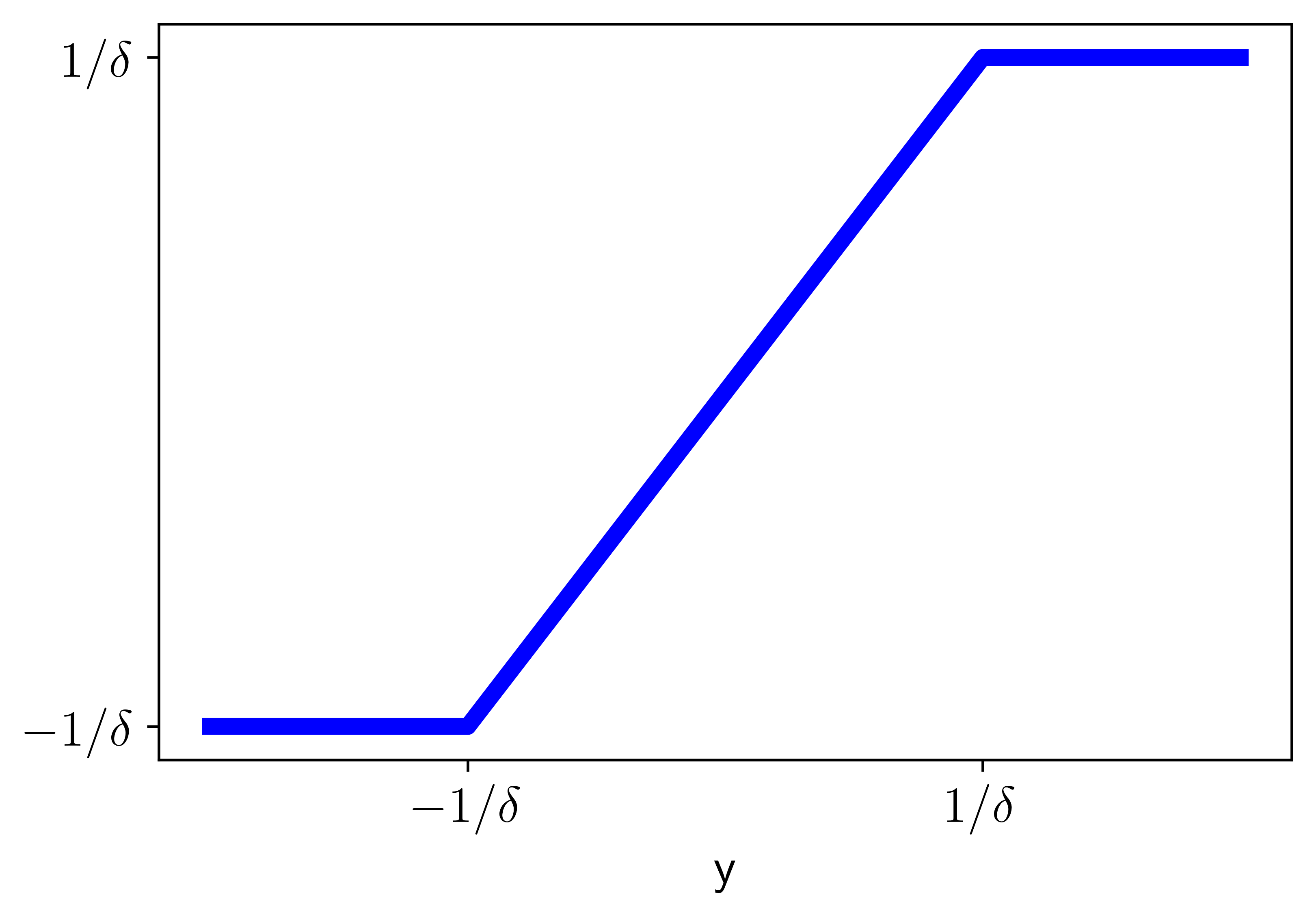}
    \caption{Cutoff $c_\delta$}
    \label{fig:cdelta}
\end{figure}

Let's start by regularizing the collision singularity in \eqref{E:initial_stochastic}. This will allow us to appeal to standard results on existence and uniqueness (see \cite{oksendal2013stochastic}).  For each $\delta>0$,  define
\begin{equation*}
c_\delta(y)
\Def \begin{cases} \tfrac1{\delta} & \text{if $y>\tfrac1{\delta}$} \\
y &\text{if $-\tfrac1{\delta}\le y\le \tfrac1{\delta}$} \\
-\tfrac1{\delta} &\text{if $y<-\tfrac1{\delta}$},  \end{cases}\end{equation*}
which is a Lipshitz-continuous truncation of the map $y\mapsto y$ at $\pm \nicefrac{1}{\delta}$ (see Figure \ref{fig:cdelta}).  Let's use standard notation
\begin{equation*} x\vee y\Def \max\{x,y\} \qquad \text{and}\qquad x\wedge y\Def \min\{x,y\} \end{equation*}
for $x$ and $y$ in $\R$. The map $(x,y)\mapsto \nicefrac{c_\delta(y)}{x^2\vee \delta^2}$ is then a bounded and Lipshitz-continuous function on $\R^2$; it is the product of two bounded Lipshitz-continuous functions.

For each $\eps>0$ and $\delta>0$,  consider the stochastic differential equation
\begin{equation*} \begin{aligned}
    d X^{\eps, \delta}_t &= Y^{\eps, \delta}_t dt\\
    d Y^{\eps, \delta}_t  &= \lb - \alpha \lb V \left(\frac{X^{\eps, \delta}_t}{d}\right)  - v_\circ +Y^{\eps, \delta}_t \rb - \beta \frac{c_\delta\left(Y^{\eps, \delta}_t\right)}{(X^{\eps,\delta}_t)^2\vee \delta^2} \rb dt  - \eps d W_t\\
    (X^{\eps, \delta}_0, Y^{\eps, \delta}_0) &= (x_\circ, y_\circ).
\end{aligned}\qquad t \ge 0 \end{equation*}
For each $\delta>0$, let's define, analogously to \eqref{E:B_dynamic}, $\BB^{(\delta)}:\R^2\to \R^2$ as
\begin{equation*}
    \BB^{(\delta)}(z) \Def \left( y,- \alpha \lb V \left(\frac{x}{d}\right) - v_\circ +y \rb - \beta \frac{c_\delta(y)}{x^2\vee \delta^2}\right); \qquad z=(x,y)\in\R^2
\end{equation*}
then $Z^{\eps,\delta}_t\Def \left(X^{\eps,\delta}_t,Y^{\eps,\delta}_t\right)$ satisfies, analogously to
\eqref{E:OV_deterministic_delta_0}, the integral equation
\begin{equation}\label{E:Sol_integ}
    Z^{\eps, \delta}_t = z_\circ + \int_0^t \BB^{(\delta)} (Z^{\eps, \delta}_s) ds - \eps \be W_t \qquad t \ge 0 
\end{equation}
for all $t>0$, where $z_\circ$ is as in \eqref{E:ZinitialDef} and $\be\Def (0,1)$.

We want to understand the dynamics of \eqref{E:Sol_integ} as a an approximation of the solution of \eqref{E:initial_stochastic}. For each $\eps>0$, we show that \eqref{E:Sol_integ} has an appropriate limit as $\delta\searrow 0$ on $\Space$ of \eqref{E:statespace}.  We do this by developing a collection of consistency results on increasing subsets of $\Space$.  This then allows us to rigorously define a solution of \eqref{E:initial_stochastic}  up to a time when the solution exits $\Space$.  Using this definition, we can mimic some of the arguments of Section \ref{S:deterministic_collision}.  We prove a probabilistic global bound using a Hamiltonian-type function.  We then construct a stochastic barrier function which will allow us to control the behavior of \eqref{E:Sol_integ} when it is near to the collision boundary $\bdy \Space$.  Appropriately combined, we can then show that collision is unlikely in \eqref{E:initial_stochastic}. 

\begin{figure}
    \centering
    \includegraphics[width=0.7\columnwidth]{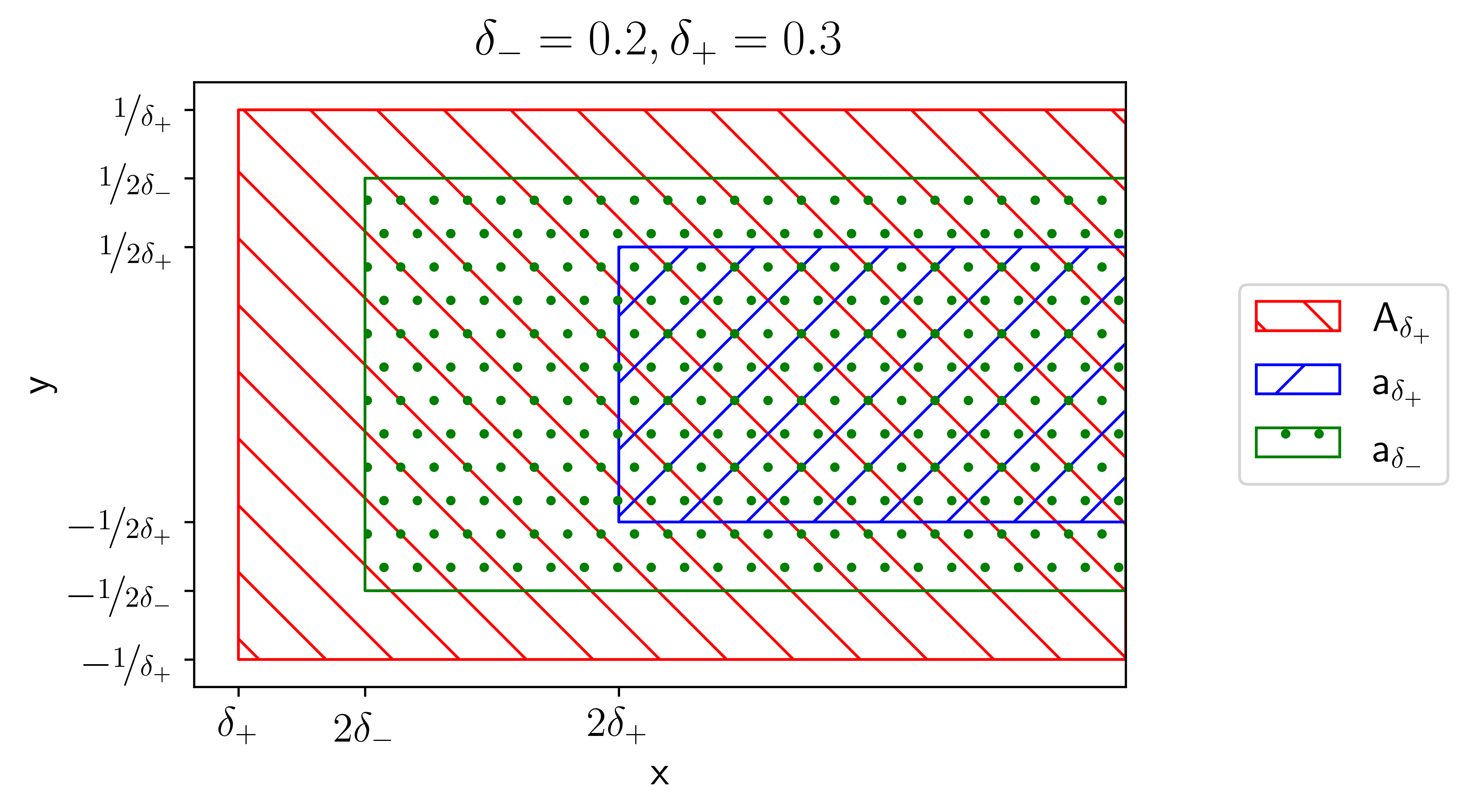}
    \caption{Sets where relaxed dynamics agree with desired ones.}
    \label{fig:agree}
\end{figure}

We carefully defined the vector field $\BB^{(\delta)}$ to agree with $\BB$ of \eqref{E:B_dynamic} on the set
\begin{equation*}\Agree_{\delta}\Def \left[\delta,\infty\right)\times \left[-\nicefrac{1}{\delta},\nicefrac{1}{\delta}\right]. \end{equation*}
Let's define the slightly smaller set
\begin{equation*} \agree_\delta \Def [2\delta,\infty)\times \left[-\nicefrac{1}{2\delta},\nicefrac{1}{2\delta}\right];\end{equation*}
see Figure \ref{fig:agree} (in fact, $\agree_\delta =\Agree_{2\delta}$, but the $\agree_{\delta}$'s and $\Agree_\delta$'s play slightly different roles, and we will thus use distinct notation).
Note that $\agree_{\delta_+}\subset \agree_{\delta_-}$ if $\delta_-<\delta_+$ (i.e., $\delta \mapsto \agree_\delta$ is decreasing in $\delta$), and 
\begin{equation*} \Space = \lim_{\delta \searrow 0}\agree_\delta = \bigcup_{\delta>0}\agree_\delta. \end{equation*}
we should be able to construct a solution of \eqref{E:initial_stochastic} by sequentially piecing together solutions of \eqref{E:Sol_integ}.  We expect, however, that a result similar to Proposition \ref{P:decreasing} should allow us to effectively restrict our calculations to a set like \eqref{E:compactset}.   We are ultimately interested in $\eps\searrow 0$ asymptotics, but want to carry out calculations on the regularized process $Z^{\delta,\eps}$.

To organize our thoughts, let's define
\begin{equation*}
    \tau^{\eps, \delta} \Def \inf \lb t\ge 0: Z^{\eps,\delta}_t\not \in \agree_\delta\rb \qquad \text{\textup{(}$\inf \emptyset\Def \infty$\textup{)}}
\end{equation*}
Since $\tau^{\eps,\delta}$ is the first time that $Z^{\eps,\delta}$ enters an open set, it is a stopping time with respect to $\{\filt_t\}_{t\ge 0}$ (which is assumed to be right-continuous).

We in fact have a \emph{consistency} result;
\begin{proposition}\label{P:consistency} for $0<\delta_-<\delta_+$, we $\BP$-a.s. have that
\begin{itemize}
    \item If $\tau^{\eps,\delta_+}<\infty$, then $\tau^{\eps,\delta_+} < \tau^{\eps,\delta_-}$
    \item If $\tau^{\eps,\delta_+}=\infty$, then $\tau^{\eps,\delta_-}=\infty$
\end{itemize}
\noindent and
\begin{equation} \label{E:consistencyclaims} \begin{aligned}
    Z^{\eps,\delta_+}_{t\wedge \tau^{\eps, \delta_+}}&=Z^{\eps,\delta_-}_{t\wedge \tau^{\eps, \delta_+}}\qquad t\ge 0 \\
    \tau^{\eps,\delta_+}&= \inf\lb t\in \left[0,\tau^{\eps,\delta_-}\right): Z^{\eps,\delta_-}_t\not\in \agree_{\delta_+}\rb. \qquad \text{\textup{(}$\inf \emptyset \Def \infty$\textup{)}}
\end{aligned}\end{equation}
\end{proposition}
\begin{proof}
Let's write
\begin{equation}\label{E:dZdifferentiable}
     Z^{\eps, \delta_+}_t-Z^{\eps, \delta_-}_t = \int_0^t \lb \BB^{(\delta_+)} (Z^{\eps, \delta_+}_s)-\BB^{(\delta_-)} (Z^{\eps, \delta_-}_s)\rb ds
\end{equation}
which shows that $Z^{\eps, \delta_+}-Z^{\eps, \delta_-}$ is differentiable in $t$.  Also note that
\begin{equation}\label{E:inclusion} \agree_{\delta_+}\subset \Agree_{\delta_+}\cap \agree_{\delta_-}\subset \Agree_{\delta_-}; \end{equation}
we will use some stopping time argument to show that $Z^{\eps,\delta_+}$ agrees with $Z^{\eps,\delta_-}$ until it leaves $\Agree_{\delta_+}\cap \agree_{\delta_-}$.  That will form the basis of the proof of the results.

Define
\begin{equation*} \bar \nu\Def \min\lb \delta_+-\delta_-,\frac{1}{\delta_-}-\frac{1}{\delta_+}\rb. \end{equation*}
Suppose that $z_+=(x_+,y_+)\in \Agree_{\delta_+}\cap \agree_{\delta_-}$ and $z_-=(x_-,y_-)$ in $\R^2$ is such that $\|z_+-z_-\|\le \bar \nu$.  Then (since $z_+\in \Agree_{\delta_+}$),
\begin{equation*} x_-\ge  x_+-\|z_+-z_-\|\ge \delta_+-\bar \nu\ge\delta_-\qquad \text{and}\qquad 
|y_-|\le |y_+|+\|z_+-z_-\|\le \frac{1}{\delta_+}+\bar \nu\le \frac{1}{\delta_-}. \end{equation*}
so $(x_-,y_-)\in \Agree_{\delta_-}$.  Then
$\BB^{(\delta_+)}(z_+)=\BB(z_+)$ and $\BB^{(\delta_-)}(z_-)=\BB(z_-)$.  Let's quantify the difference between $\BB(z_-)$ and $\BB(z_+)$ by defining
\begin{equation}\label{E:consistencyK} \KK_{\eqref{E:consistencyK}}\Def \sup_{\substack{z_+\in \Agree_{\delta_+}\cap \agree_{\delta_-} \\ z_-\in \Agree_{\delta_-} \\ z_-\not = z_+}}\frac{\left\|\BB(z_-)-\BB(z_+)\right\|}{\left\|z_--z_+\right\|} \end{equation}
(which is bounded from above by the Lipshitz coefficient of $\BB\big|_{\Agree_{\delta_-}}$).

Fix now $\nu\in (0,\bar \nu)$, $L>0$, and $\vsig>0$ (the parameters $\eps$, $\delta_+$ and $\delta_-$ are fixed in this proof), and define
\begin{align*} \sigma_1 &= \inf \lb t\ge 0: Z^{\eps,\delta_+}_t\not \in \Agree_{\delta_+}\cap \agree_{\delta_-}\rb\\
\sigma_{2,\nu} &= \inf \lb t\ge 0: \left\|Z^{\eps,\delta_+}_t-Z^{\eps,\delta_-}_t\right\|>\nu\rb\wedge \sigma_1\\
E^{\vsig}_t&\Def \lb \left\|Z^{\eps,\delta_+}_t-Z^{\eps,\delta_-}_t\right\|^2+\vsig^2\rb^{1/2}e^{-\KK_{\eqref{E:consistencyK}} t} \qquad t\ge 0 \end{align*}

Differentiating $E^{\vsig}$ with respect to time (recall \eqref{E:dZdifferentiable}), we have that
\begin{align*} \dot E^{\vsig}_t &=\lb \frac{\la Z^{\eps,\delta_+}_t-Z^{\eps,\delta_-}_t,\BB^{(\delta_+)}(Z^{\eps,\delta_+}_s)-\BB^{(\delta_-)}(Z^{\eps,\delta_-}_t)\ra}{\lb \left\|Z^{\eps,\delta_+}_t-Z^{\eps,\delta_-}_t\right\|^2+\vsig^2\rb^{1/2}} \right.\\
&\qquad \left. -\KK_{\eqref{E:consistencyK}}\lb \left\|Z^{\eps,\delta_+}_t-Z^{\eps,\delta_-}_t\right\|^2+\vsig^2\rb^{1/2}\rb e^{-\KK_{\eqref{E:consistencyK}}t} \end{align*}

If $0<t<\sigma_{2,\nu}$,
\begin{align*} \dot E^{\vsig}_t&\le \lb \frac{\KK_{\eqref{E:consistencyK}}\left\| Z^{\eps,\delta_+}_t-Z^{\eps,\delta_-}_t\right\|^2}{\lb \left\|Z^{\eps,\delta_+}_t-Z^{\eps,\delta_-}_t\right\|^2+\vsig^2\rb^{1/2}}\right.\\
&\qquad \left. -\KK_{\eqref{E:consistencyK}}\lb \left\|Z^{\eps,\delta_+}_t-Z^{\eps,\delta_-}_t\right\|^2+\vsig^2\rb^{1/2}\rb e^{-\KK_{\eqref{E:consistencyK}}t}\le 0 \end{align*}
and thus
\begin{equation*} \left\| Z^{\eps,\delta_+}_{\sigma_{2,\nu}\wedge L}-Z^{\eps,\delta_-}_{\sigma_{2,\nu}\wedge L}\right\|\exp\left[-\KK_{\eqref{E:consistencyK}}\lb \sigma_{2,\nu}\wedge L\rb \right]\le E^\vsig_{\sigma_{2,\nu}}\le E^\vsig_0 = \vsig. \end{equation*}
Rearranging, we get that
\begin{equation*} \left\| Z^{\eps,\delta_+}_{\sigma_{2,\nu}\wedge L}-Z^{\eps,\delta_-}_{\sigma_{2,\nu}\wedge L}\right\|\le \vsig \exp\left[\KK_{\eqref{E:consistencyK}}\lb \sigma_{2,\nu}\wedge L\rb \right]. \end{equation*}
Letting $\vsig\searrow 0$, we get that
\begin{equation*} \left\| Z^{\eps,\delta_+}_{\sigma_{2,\nu}\wedge L}-Z^{\eps,\delta_-}_{\sigma_{2,\nu}\wedge L}\right\|=0. \end{equation*}
Finally letting $L\nearrow \infty$, we get that
\begin{equation*} \left\| Z^{\eps,\delta_+}_{\sigma_{2,\nu}}-Z^{\eps,\delta_-}_{\sigma_{2,\nu}}\right\|=0. \end{equation*}
on the set where $\sigma_{2,\nu}<\infty$.  This is impossible if $\sigma_{2,\nu}<\sigma_1$, so $\sigma_{2,\nu}=\sigma_1$ ($\BP$-a.s.), so 
\begin{equation*} \sup_{0\le t< \sigma_1}\left\| Z^{\eps,\delta_+}_t-Z^{\eps,\delta_-}_t\right\|\le \nu. \end{equation*}
Letting $\nu\searrow 0$, we have that
\begin{equation}\label{E:consistencyequals} \sup_{0\le t< \sigma_1}\left\| Z^{\eps,\delta_+}_t-Z^{\eps,\delta_-}_t\right\|=0. \end{equation}

We next claim that 
\begin{equation} \label{E:strictagree} \tau^{\eps,\delta_+}<\sigma_1\le \tau^{\eps,\delta_-}. \end{equation}
If $t<\tau^{\eps,\delta_+}$, then $Z^{\eps,\delta_+}_t\in \agree_{\delta_+}$.  In light of the first inclusion in \eqref{E:inclusion}, we have that $\sigma_1\ge \tau^{\eps,\delta_+}$.   If $\tau^{\eps,\delta_+}<\infty$, then $Z^{\eps,\delta_+}_{\tau^{\eps,\delta_+}}=Z^{\eps,\delta_+}_{\tau^{\eps,\delta_+}-}\in \agree_{\delta_+}$ (by continuity of $Z^{\eps,\delta_+}$ and the fact that $\agree_{\delta_+}$ is closed).  In fact, $\agree_{\delta_+}$ is contained in the interior of 
$\Agree_{\delta_+}\cap \agree_{\delta_-}$, so the continuity of $Z^{\eps,\delta_+}$ implies that there is a $\nu>0$ such that $Z^{\eps,\delta_+}_t\in \Agree_{\delta_+}\cap \agree_{\delta_-}$ for $t\in [\tau^{\eps,\delta_+},\tau^{\eps,\delta_+}+\nu)$; thus $\sigma_1 \ge \tau^{\eps,\delta_+}+\nu$.  The left-hand inequality of \eqref{E:strictagree} follows.  
If $t<\sigma_1$, then $Z^{\eps,\delta_-}_t=Z^{\eps,\delta_+}_t\in \Agree_{\delta_+}\cap \agree_{\delta_-}\subset \agree_{\delta_-}$ (in light of \eqref{E:consistencyequals}); so the right-hand claim of \eqref{E:strictagree} follows. The chain \eqref{E:strictagree} of inequalities directly proves the first two claim of Proposition \ref{P:consistency}, and, when combined with \eqref{E:consistencyequals}, also gives us the first claim of \eqref{E:consistencyclaims}.

Define
\begin{equation*} \hat \tau\Def \inf\lb t\in \left[0,\tau^{\eps,\delta_-}\right): Z^{\eps,\delta_-}_t\not\in \agree_{\delta_+}\rb\end{equation*}
with the right-hand side being $\infty$ upon taking $\inf \emptyset$ (the standard convention).
If $t<\tau^{\eps,\delta_+}$ (and thus in turn $t<\sigma_1$), $Z^{\eps,\delta_-}_t=Z^{\eps,\delta_+}_t\in \agree_{\delta_+}$ (from the first claim of \eqref{E:consistencyclaims}); thus $\hat \tau\ge \tau^{\eps,\delta_+}$.  From \eqref{E:strictagree}, we conversely have that
\begin{equation*} \tau^{\eps,\delta_+}
=\inf\lb t\in [0,\sigma_1): Z^{\eps,\delta_+}_t\not\in \agree_{\delta_+}\rb 
= \inf\lb t\in [0,\sigma_1): Z^{\eps,\delta_-}_t\not\in \agree_{\delta_+}\rb \ge \hat \tau \end{equation*}
with these inequalities trivially holding with the convention that $\inf \emptyset\Def \infty$.
This gives us the last claim of \eqref{E:consistencyclaims}.
\end{proof}

Let's now define
\begin{equation}\label{E:tau_eps}
    \tau^\eps \Def \sup_{m\ge 1}\tau^{\eps,\nicefrac{1}{m}};
\end{equation}
since $\tau^\eps$ is a supremum of a countable collection of stopping times, it too is a stopping time.
In light of the first claims of Proposition \ref{P:consistency}, $\tau^\eps=\lim_{m\nearrow \infty}\tau^{\eps,\nicefrac{1}{m}}$, $\BP$-a.s.
Let's next piece together the $Z^{\eps,\delta}$'s on $[0,\tau^\eps)$, being careful to do so in a way which will preserve various properties \cite{ethier2009markov}.  Let's add a cemetery state to formalize what happens after $\tau^\eps$.  Namely, fix a point $\star$ not in $\Space$, and define
\begin{equation*} \Space^\star \Def \Space\cup \{\star\}. \end{equation*}
and endow $\Space$ with the standard topology of one-point compactifications; this will allow us to rigorously frame our existence and uniqueness results on $\Space$ of \eqref{E:Sol_integ}. 
For each positive integer $m$, let's then define
\begin{equation*} \hat Z^{\eps,m}_t\Def \begin{cases} Z^{\eps,\nicefrac{1}{m}}_t &\text{if $t<\tau^{\eps,\nicefrac{1}{m}}$} \\
\star &\text{if $t\ge \tau^{\eps,\nicefrac{1}{m}}$} \end{cases} \end{equation*}
allowing us to separate killing from the singularity at $\bdy \Space$.
\begin{theorem}\label{T:Limit}
Fix $\eps>0$.
For each $t\ge 0$, $\hat Z^\eps_t\Def \lim_{m\nearrow \infty}\hat Z^{\eps,m}_t$ is $\BP$-a.s. well-defined \textup{(}in the topology of $\Space^\star$\textup{)}.  For each $\delta>0$,
we $\BP$-a.s. have that
\begin{equation}\label{E:LimitClaims}\begin{aligned} \hat Z^\eps_{t\wedge \tau^{\eps,\delta}} &= Z^{\eps,\delta}_{t\wedge \tau^{\eps,\delta}}. \qquad t\ge 0 \\
\hat Z^\eps_{t\wedge \tau^{\eps,\delta}} &= z_\circ + \int_{s=0}^{t\wedge \tau^{\eps,\delta}}\BB(\hat Z^\eps_s)ds - \eps \be W_{t\wedge \tau^{\eps,\delta}}\\
\tau^{\eps,\delta}&= \inf\lb t\in [0,\tau^\eps): \hat Z^\eps_t\not\in \agree_\delta\rb. \qquad \text{ \textup{(}$\inf \emptyset\Def\infty$\textup{)}}\end{aligned}\end{equation}
\end{theorem}
\noindent In particular, $\hat Z^\eps$ is the solution of \eqref{E:initial_stochastic} on $[0, \tau^\eps)$.
\begin{proof}[Proof of Theorem \ref{T:Limit}]
If $t<\tau^\eps$, then $t<\tau^{\eps,\nicefrac{1}{m}}$ for some positive integer $m$ (the first claim of Proposition \ref{P:consistency} being in this case a strict inequality).  For each positive integer $m'> m$, $t<\tau^{\eps,\nicefrac{1}{m}}<\tau^{\eps,\nicefrac{1}{m'}}$ $\BP$-a-s., so
\begin{equation*} \hat Z^{\eps,m'}_t = Z^{\eps,\nicefrac{1}{m'}}_t=Z^{\eps,\nicefrac{1}{m}}_t  \end{equation*}
$\BP$-a.s. (use the definition of $\hat Z^{\eps,m'}$ and then the first claim of \eqref{E:consistencyclaims}). Thus $\hat Z^\eps_t$ is well-defined if $t<\tau^\eps$.  If $t\ge \tau^\eps$, then $t>\tau^{\eps,\nicefrac{1}{m}}$ $\BP$-a.s. for all positive integers $m$ (the inequality in the first claim of Proposition \ref{P:consistency} being strict), so $\hat Z^{\eps,m}_t=\star$ $\BP$-a.s. for all positive integers $m$.  Thus $\hat Z^\eps_t$ is well-defined if $t\ge \tau^\eps$.

For $\delta>0$ and then any positive integer $m$ with  $m>\nicefrac{1}{\delta}$ (so that $\nicefrac{1}{m}<\delta$ and hence $\tau^{\eps,\delta}<\tau^{\eps,\nicefrac{1}{m}}$ $\BP$-a.s. by the first claim of Proposition \ref{P:consistency}, 
we then have that
\begin{equation}\label{E:BBB} \hat Z^{\eps,m}_{t\wedge \tau^{\eps,\delta}} = Z^{\eps,\nicefrac{1}{m}}_{t\wedge \tau^{\eps,\delta}}=Z^{\eps,\delta}_{t\wedge \tau^{\eps,\delta}} \end{equation}
for each $t\ge 0$, $\BP$-a.s. (using the first claim of \eqref{E:consistencyclaims}).  Taking limits in $m$, we get the first claim of \eqref{E:LimitClaims}.  From \eqref{E:BBB} and the dynamics of $Z^{\eps,\delta}$, we have that
\begin{align*} \hat Z^\eps_{t\wedge \tau^{\eps,\delta}} 
&= Z^{\eps,\delta}_{t\wedge \tau^{\eps,\delta}}
= z_\circ + \int_{s=0}^{t\wedge \tau^{\eps,\delta}}\BB^{(\delta)}(Z^{\eps,\delta}_s)ds - \eps \be W_{t\wedge \tau^{\eps,\delta}}\\
&= z_\circ + \int_{s=0}^{t\wedge \tau^{\eps,\delta}}\BB(Z^{\eps,\delta}_s)ds - \eps \be W_{t\wedge \tau^{\eps,\delta}}
= z_\circ + \int_{s=0}^{t\wedge \tau^{\eps,\delta}}\BB(\hat Z^\eps_s)ds - \eps \be W_{t\wedge \tau^{\eps,\delta}}\end{align*}
which is the second claim of \eqref{E:LimitClaims}.

Define
\begin{equation*} \hat \tau\Def \inf\lb t\in [0,\tau^\eps): \hat Z^\eps_t\not\in \agree_\delta\rb  \end{equation*}
with the right-hand side being $\infty$ upon taking $\inf \emptyset$ (the standard convention).
If $t<\tau^{\eps,\delta}$, then the first claim of \eqref{E:LimitClaims} implies that $\hat Z^\eps_t=Z^{\eps,\delta}_t\in \agree_\delta$, implying that $\hat \tau\ge \tau^{\eps,\delta}$. Let's again fix an integer $m>\nicefrac{1}{\delta}$ (so that $\delta>\nicefrac{1}{m}$).  From the last claim of \eqref{E:consistencyclaims} and the first claim of \eqref{E:LimitClaims}), we have that
\begin{equation*} \tau^{\eps,\delta}
=\inf\lb t\in \left[0,\tau^{\eps,\nicefrac{1}{m}}\right): Z^{\eps,\nicefrac{1}{m}}_t\not\in \agree_{\delta}\rb
=\inf\lb t\in \left[0,\tau^{\eps,\nicefrac{1}{m}}\right): \hat Z^\eps_t\not\in \agree_{\delta}\rb
\ge \hat \tau, \end{equation*}
with these inequalities trivially holding with the convention that $\inf \emptyset \Def \infty$.
This gives us the last claim of \eqref{E:LimitClaims}.
\end{proof}
Since the limit $\hat Z^\eps_t$ is well-defined $\BP$-a.s. for each $t$, $\hat Z^\eps_t$ is adapted.  For $t\in [0,\tau^\eps)$, let's write $\hat Z^\eps_t=(X^\eps_t,Y^\eps_t)$; the pair $(X^\eps,Y^\eps)$ satisfies \eqref{E:initial_stochastic} up to (but not including) $\tau^\eps$.

To proceed, let's prove a stochastic complement to Proposition \ref{P:decreasing} (and recall  \eqref{E:hcircdef}).  For $\eps>0$, define
\begin{equation}\label{E:stopTime_bound}
    \tau^{H,\eps} \Def \inf \lb t \in [0,\tau^\eps): H(X^\eps_t.Y^\eps_t)>h_\circ+1 \rb \qquad \text{(with $\inf \emptyset \Def \infty$)}
\end{equation}
Since the $\tau^{\eps,\nicefrac{1}{m}}$'s are strictly increasing in $m$, 
\begin{equation}\label{E:aaa} \lb \tau^{H,\eps}<t\rb= \bigcup_{\substack{s\in \Q \\ s<t}}\bigcup_{m=1}^\infty \left(\lb H(X^\eps_s.Y^\eps_s)>h_\circ+1\rb \cap \{\tau^{\eps,\nicefrac{1}{m}}>s\}\right)\in \filt_t \end{equation}
(with $\Q$ being the set of rational numbers)
for each $t>0$, and hence (using the assumption that the filtration $\{\filt\}_{t\ge 0}$ is right-continuous)
\begin{equation}\label{E:aab} \lb \tau^{H,\eps}\le t\rb = \bigcap_{\substack{t'>t \\ t'\in \Q}}\lb \tau^{H,\eps}<t'\rb \in \filt_{t+}=\filt_t \end{equation}
so $\tau^{H,\eps}$ is indeed a stopping time with respect to the filtration $\{\filt_t\}_{t\ge 0}$.
In light of \eqref{E:levelsetscompact}, we thus have that
\begin{equation}\label{E:YisBounded} |Y^\eps_t|\le \bar y \end{equation}
for $0\le t<\tau^{H,\eps}\wedge \tau^\eps$.
\begin{theorem} \label{T:Hbounded}
We have that
\begin{equation*}
    \lim_{\substack{\eps \to 0 \\ \eps \sqrt{L} \to 0}} \BP \lb \tau^{H,\eps} <L \rb = 0.
\end{equation*}
\end{theorem}
\noindent In other words, as $\eps \searrow 0$ and $L$ potentially becomes large, but such that $L = o \left(\eps^{-2}\right)$, where $o(\cdot)$ is little o notation, then $\BP \lb \tau^{H,\eps} <L \rb$ becomes small.
\begin{proof}[Proof of Theorem \ref{T:Hbounded}]
Using the second claim of \eqref{E:LimitClaims}, let's apply Ito's formula to the Hamiltonian function $H$.  Fix a positive integer $m$; we get
(similarly to \eqref{E:energyLevel})
\begin{equation}\label{E:Ito}
     H\left(X^\eps_{t \wedge \tau^{\eps,\nicefrac{1}{m}}}, Y^\eps_{t \wedge \tau^{\eps,\nicefrac{1}{m}}}\right) = h_\circ - \int_{s=0}^{t \wedge \tau^{\eps,\nicefrac{1}{m}}}\lb  \alpha  (Y^\eps_s)^2  + \frac{\beta Y^\eps_s c_\delta(Y^\eps_s)}{(X^\eps_s)^2\vee \left(\nicefrac{1}{m}\right)^2} \rb ds - \eps \int_{s=0}^{t \wedge \tau^{\eps,\nicefrac{1}{m}}}  Y^\eps_s dW_s + \tfrac 12 \eps^2 (t \wedge \tau^{\eps, \nicefrac{1}{m}}), 
\end{equation}
for $t\ge 0$.  If $\tau^{H,\eps}<L\wedge \tau^{\eps, \nicefrac{1}{m}}$, then $H\left(X^\eps_{\tau^{H,\eps}\wedge L\wedge \tau^{\eps,\nicefrac{1}{m}}},Y^\eps_{\tau^{H,\eps}\wedge L\wedge \tau^{\eps,\nicefrac{1}{m}}}\right)\ge h_\circ+1$.
Using \eqref{E:Ito}, we can write
\begin{align*}
    H\left(X^\eps_{\tau^{H,\eps}\wedge L\wedge \tau^{\eps,\nicefrac{1}{m}}},X^\eps_{\tau^{H,\eps}\wedge L\wedge \tau^{\eps,\nicefrac{1}{m}}}\right) &\le  h_\circ - \eps \int_{s=0}^{\tau^{H,\eps}\wedge L\wedge \tau^{\eps,\nicefrac{1}{m}}} Y^\eps_s dW_s + \frac 12 \eps^2 (\tau^{H,\eps}\wedge L\wedge \tau^{\eps,\nicefrac{1}{m}}) \\
    & \le h_\circ + \eps \left|\int_{s=0}^{\tau^{H,\eps}\wedge L\wedge \tau^{\eps,\nicefrac{1}{m}}} Y^\eps_s dW_s \right|+ \frac 12 \eps^2 (\tau^{H,\eps}\wedge L\wedge \tau^{\eps,\nicefrac{1}{m}})\\
    & \le h_\circ + \eps \left|\int_{s=0}^{\tau^{H,\eps}\wedge L\wedge \tau^{\eps,\nicefrac{1}{m}}} Y^\eps_s dW_s \right|+ \frac 12 \eps^2 L;
\end{align*}
we have used here the fact that $(Y^\eps_s)^2$ and $Y^\eps_s c_\delta(Y^\eps_s)$ in \eqref{E:Ito} are nonnegative.
Combining things and assuming that $\eps^2 L\le 1$, we have that
\begin{equation*}
    \lb \tau^{H,\eps}<L\wedge \tau^{\eps, \nicefrac{1}{m}} \rb  \subset  \lb h_\circ + \eps \left|\int_0^{\tau^{H,\eps}\wedge L\wedge \tau^{\eps,\nicefrac{1}{m}}} Y^\eps_s dW_s \right|+ \frac 12  \ge 1 + h_\circ \rb
    \subset \lb  \eps \left|\int_0^{\tau^{H,\eps}\wedge L\wedge \tau^{\eps,\nicefrac{1}{m}}} Y^\eps_s dW_s \right| \ge \frac 12\rb.
\end{equation*}
Using \eqref{E:YisBounded}, we can further write
\begin{align*}
    \BP \lb \tau^{H,\eps}<L\wedge \tau^{\eps, \nicefrac{1}{m}} \rb
    &\le  \BP \lb \eps  \left|\int_0^{\tau^{H,\eps}\wedge L\wedge \tau^{\eps,\nicefrac{1}{m}}} Y^\eps_s dW_s \right| \ge \frac 12\rb
    \le 4 \eps^2 \BE \left[ \left(\int_{s=0}^{\tau^{H,\eps}\wedge L\wedge \tau^{\eps,\nicefrac{1}{m}}} Y^\eps_s dW_s \right)^2 \right]\\
    & = 4 \eps^2 \BE \left[\int_0^{\tau^{H,\eps}\wedge L\wedge \tau^{\eps,\nicefrac{1}{m}}} \left(Y^\eps_s\right)^2 ds  \right]
    \le 4 \eps^2 \bar y^2 L,
\end{align*}
where equality on the right follows from Ito's isometry.  Letting $m \to \infty$, the claim follows.
\end{proof}

Our main result shows that collisions are unlikely, with the same type of asymptotics as Theorem \ref{T:Hbounded} (and indeed the proof will use Theorem \ref{T:Hbounded}).
\begin{theorem}\label{T:main} Let $\tau^\eps$ be as of \eqref{E:tau_eps}. We have that
\begin{equation*} 
 \lim_{\substack{\eps \to 0 \\ \eps \sqrt{L} \to 0}} \BP \lb \tau^\eps <L \rb = 0.\end{equation*}
 \end{theorem}
\noindent As with Theorem \ref{T:Hbounded}, we want to do so by applying the tools of stochastic analysis to the middle claim of \eqref{E:LimitClaims}.  

Our proof has several key components. Firstly, \eqref{E:YisBounded} tells us that we can essentially localize our calculations to $\R\times [-\bar y,\bar y]$.  Secondly, the ideas of \eqref{E:local_parametrization} and \eqref{E:localpara} suggest that we compare $(X^\eps,Y^\eps)$ to a (vertical) \emph{reference manifold} parametrized by $Y^\eps$; the ODE \eqref{E:phi_dynamic_deterministic} will serve as our guide in constructing it.  We will then carefully create a (smooth) \emph{danger} function which quantifies proximity to collision by comparing it to the reference manifold.  Mathematically, the danger function will be nonincreasing along curves of the dominant (deterministic) part of \eqref{E:initial_stochastic}; this will allow us to bound the likelihood of collision from above.

Let's start by setting
\begin{equation} \label{E:xdaggerdef} x_\dagger\Def \min\lb x_\circ,dV^{-1}\left(\nicefrac{v_\circ}{2}\right)\rb; \end{equation}
this is similar to \eqref{E:xminusdef}, but allows for an extra margin around $x_\infty$.
Consider the ODE
\begin{equation} \label{E:stochasticbarrier}\begin{aligned} \phi_\dagger'(y)&= -\frac{\phi_\dagger^2(y)}{\alpha \phi_\dagger^2(y)+\beta}\qquad y\in \R \\
\phi_\dagger(-\bar y)&=x_\dagger. \end{aligned}\end{equation}
(which is the analogue of \eqref{E:phi_dynamic_deterministic}).
Let's also define
\begin{equation} \label{E:uphidef} \uphi\Def \left(\frac{1}{x_\dagger}+\frac{2\bar y}{\beta}\right)^{-1}. \end{equation}
\begin{lemma}\label{L:stochasticbarrier} The ODE \eqref{E:stochasticbarrier} is well-defined on $\R$.  The solution $\phi_\dagger$ is decreasing. Furthermore
\begin{equation} \label{E:philowerbound} \inf_{|y|\le \bar y}\phi_\dagger(y)\ge \uphi \end{equation}
\textup{(}with $\uphi$ as in \eqref{E:uphidef}\textup{)}.

There is a $\KK_{\eqref{E:barrierbounds}}>0$ such that
\begin{equation}\label{E:barrierbounds} \left|\phi_\dagger(y)\right|\le \KK_{\eqref{E:barrierbounds}} \qquad \text{and}\qquad \left|\phi_\dagger'(y)\right|\le \KK_{\eqref{E:barrierbounds}} \qquad \text{and}\qquad \left|\phi_\dagger^{\prime \prime}(y)\right|\le \KK_{\eqref{E:barrierbounds}} \end{equation}
for all $y\in [-\bar y,\bar y]$.
\end{lemma}
\noindent A lower bound on $\phi_\dagger$ in fact naturally follows by uniqueness and the fact that the only fixed point of \eqref{E:stochasticbarrier} is at $0$; the statement of Lemma \ref{L:stochasticbarrier} gives some more precise bounds.
\begin{proof}[Proof of Lemma \ref{L:stochasticbarrier}] Define
\begin{equation}\label{E:frhsdef} f(y)\Def -\frac{y^2}{\alpha y^2 + \beta} \qquad y\in \R \end{equation}
so that the ODE \eqref{E:stochasticbarrier} is
\begin{equation}\label{E:barrierODE} \phi_\dagger'(y)=f(\phi_\dagger(y)). \end{equation}
Note that
\begin{equation}\label{E:fbound} \left|f(y)\right|\le \frac{y^2}{\alpha y^2}=\frac{1}{\alpha} \end{equation}
for $y\in \R$ so $f$ is in fact bounded.
Next note that
\begin{equation*} f'(y) = -\frac{2y}{\alpha y^2+\beta}+\frac{2\alpha y^3}{\left(\alpha y^2 + \beta\right)^2} \qquad y\in \R^2 \end{equation*}
Young's inequality implies that
\begin{align*} |y|&=\frac{1}{\sqrt{\alpha \beta}}\left|(y\sqrt{\alpha})(\sqrt{\beta})\right|
\le \frac{1}{2\sqrt{\alpha \beta}}\lb \alpha y^2 +\beta\rb \\
|y|^{3/2} &= \frac{1}{\alpha^{3/4}(3 \beta)^{1/4}}\left|(y\sqrt{\alpha})^{3/2}(3\beta)^{1/4}\right|
\le \frac{1}{\alpha^{3/4}(3 \beta)^{1/4}}\lb \frac34 (y\sqrt{\alpha})^{3/2\times 4/3} +\frac14 (3\beta)^{1/4\times 4}\rb \\
&= \frac14\left(\frac{3^3}{\alpha^3\beta}\right)^{1/4}\lb \alpha y^2+ \beta\rb
\end{align*}
for all $y\in \R$.
Thus
\begin{equation}\label{E:fprimebound} |f'(y)| \le  \frac{2|y|}{\alpha y^2+\beta}+\frac{(2\alpha)\lb|y|^{3/2}\rb^2}{\left(\alpha y^2 + \beta\right)^2}\le \frac{1}{\sqrt{\alpha \beta}}\frac{\alpha y^2 + \beta}{\alpha y^2 + \beta} + \frac{\alpha}{8}\left(\frac{3^3}{\alpha^3\beta}\right)^{1/2}\frac{(\alpha y^2+\beta)^2}{(\alpha y^2 + \beta)^2} \le \frac{1}{\sqrt{\alpha \beta}}+\frac{\alpha}{8}\left(\frac{3^3}{\alpha^3\beta}\right)^{1/2} \end{equation}
for all $y\in \R$, implying that $f$ is Lipshitz.  We now have that \eqref{E:stochasticbarrier} is well-defined on all of $\R$.  Since $f$ of \eqref{E:frhsdef} is negative, the solution $\phi$ of \eqref{E:barrierODE} is decreasing.

Define
\begin{equation*} y_+\Def \inf\lb y \in [-\bar y,\infty): \phi_\dagger(y)=0\rb \qquad \text{($\inf\emptyset\Def \infty$)} \end{equation*}
For $y\in (-\bar y,y_+)$, 
\begin{equation*} \frac{d}{dy}\phi_\dagger^{-1}(y)=-\frac{\dot \phi_\dagger(y)}{\phi_\dagger^2(y)} = \frac{1}{\alpha \phi_\dagger^2(y)+\beta}\le \frac{1}{\beta}. \end{equation*}
 For $y\in [-\bar y,y_+)$, we thus have that
\begin{equation*} \frac{1}{\phi_\dagger(y)}-\frac{1}{\phi_\dagger(-\bar y)}\le \frac{1}{\beta}\left(y+\bar y\right) \end{equation*}
(compare with the proof of Lemma \ref{L:collision_avoid}) which implies that
\begin{equation}\label{E:lowerbarrier}\phi_\dagger(y)>\left(\frac{1}{x_\dagger}+\frac{y+\bar y}{\beta}\right)^{-1} \end{equation}
for $y\in [-\bar y,y_+)$.  
If $y_+<\infty$, we can take $y\nearrow y_+$ in \eqref{E:lowerbarrier} and arrive at the contradiction
\begin{equation*} 0=\lim_{y\nearrow y_+}\phi_\dagger(y)\ge \left(\frac{1}{x_\dagger}+\frac{y_++\bar y}{\beta}\right)^{-1}. \end{equation*}
Thus $y_+=\infty$ and we consequently have \eqref{E:lowerbarrier} for all $y\ge -\bar y$; this implies \eqref{E:philowerbound}.

To see the second bound of \eqref{E:barrierbounds}, combine
\eqref{E:barrierODE} and \eqref{E:fbound}.  To see the third bound of \eqref{E:barrierbounds}, differentiate \eqref{E:barrierODE} to get
\begin{equation*} \phi_\dagger^{\prime \prime}(y) = f'(\phi_\dagger(y))\phi_\dagger'(y)= f'(\phi_\dagger(y))f(\phi_\dagger(y)) \end{equation*}
and then combine \eqref{E:fbound} and \eqref{E:fprimebound}.  The first bound of \eqref{E:barrierbounds} comes from the second bound and the boundary condition of \eqref{E:stochasticbarrier}.
\end{proof}

Let's now flatten $\phi_\dagger$ to the right of zero, with a transition region which depends on various system parameters.  Define
\begin{equation}\label{E:varpidef}
    \varpi_\dagger \Def \frac{\nicefrac{\alpha v_\circ}{2}}{\alpha + \nicefrac{4\beta}{\uphi^2}}.
\end{equation}
Let $\varrho:\R\to [0,1]$ in $C^\infty$ be nonincreasing function and such that $\varrho:(-\infty,0]\mapsto 1$ and $\varrho:[1,\infty)\mapsto 0$.
Define
\begin{equation}\label{E:barrierdef} \phi(y)\Def \varrho\left(\frac{y}{\varpi_\dagger}\right)\phi_\dagger(y) + \lb 1- \varrho\left(\frac{y}{\varpi_\dagger}\right)\rb \phi_\dagger(\varpi_\dagger) = \varrho\left(\frac{y}{\varpi_\dagger}\right)\lb \phi_\dagger(y)-\phi_\dagger(\varpi_\dagger)\rb + \phi_\dagger(\varpi_\dagger) \end{equation}
\begin{lemma}\label{L:phiprops} The function $\phi$ is nonincreasing and $\inf_{|y|\le \bar y}\phi(y)\ge \uphi$, where $\uphi$ is defined in \eqref{E:uphidef}.  We have that
\begin{equation}\label{E:phibarder} \phi'(y) =\begin{cases} -\frac{1}{\alpha + \beta/\phi^2(y)} &\text{if $y\in [-\bar y,0]$} \\
 0 &\text{if $y\in [\varpi_\dagger,\bar y]$.}  \end{cases}\end{equation}

There is a $\KK_{\eqref{E:phibds}}>0$ such that
\begin{equation}\label{E:phibds} \left|\phi(y)\right|\le \KK_{\eqref{E:phibds}} \qquad \text{and}\qquad \left|\phi'(y)\right|\le \KK_{\eqref{E:phibds}} \qquad \text{and}\qquad \left|\phi^{\prime \prime}(y)\right|\le \KK_{\eqref{E:phibds}} \end{equation}
for $y\in [-\bar y,\bar y]$.
\end{lemma}
\noindent This result still holds in the case that $\varpi_\dagger>\bar y$.
\begin{proof}[Proof of Lemma \ref{L:phiprops}]
Since $\phi_\dagger$ is nonincreasing, $y\mapsto \phi_\dagger(y)-\phi_\dagger(\varpi_\dagger)$ is nonnegative on $[-\bar y,\varpi_\dagger]$.  Since $\varrho$ is also nonnegative, the second representation of \eqref{E:barrierdef} implies that $\phi(y)\ge \phi_\dagger(\varpi_\dagger)$ for $y\in [-\bar y,\varpi_\dagger]$. Since $\phi(y)$ is constant on $[\varpi_\dagger,\bar y]$, the second claim follows from \eqref{E:lowerbarrier}.

The first bound of \eqref{E:phibds} follows from the corresponding bound of \eqref{E:barrierbounds}
Taking one then two derivatives of $\phi$, we have
\begin{equation} \label{E:barrierderivatives}\begin{aligned} \phi'(y) &=\varrho\left(\frac{y}{\varpi_\dagger}\right)\phi_\dagger'(y)+ \frac{1}{\varpi_\dagger}\varrho'\left(\frac{y}{\varpi_\dagger}\right)\lb \phi_\dagger(y)-\phi_\dagger(\varpi_\dagger)\rb \\
\phi^{\prime \prime}(y) &= \varrho\left(\frac{y}{\varpi_\dagger}\right)\phi_\dagger^{\prime \prime}(y)+ \frac{2}{\varpi_\dagger}\varrho'\left(\frac{y}{\varpi_\dagger}\right)\phi'_\dagger(y)+\frac{1}{\varpi^2_\dagger}\varrho^{\prime \prime}\left(\frac{y}{\varpi_\dagger}\right)\lb \phi_\dagger(y)-\phi_\dagger(\varpi_\dagger)\rb\end{aligned}\end{equation}

The claimed derivatives of \eqref{E:phibarder} follow from the first equality of \eqref{E:barrierderivatives} which also implies that $\phi$ is nonincreasing. The bounds on the first two derivatives of $\phi$ follow from \eqref{E:barrierderivatives} and \eqref{E:barrierbounds}.
\end{proof}

\begin{figure}
    \centering
    \includegraphics[width=0.7\columnwidth]{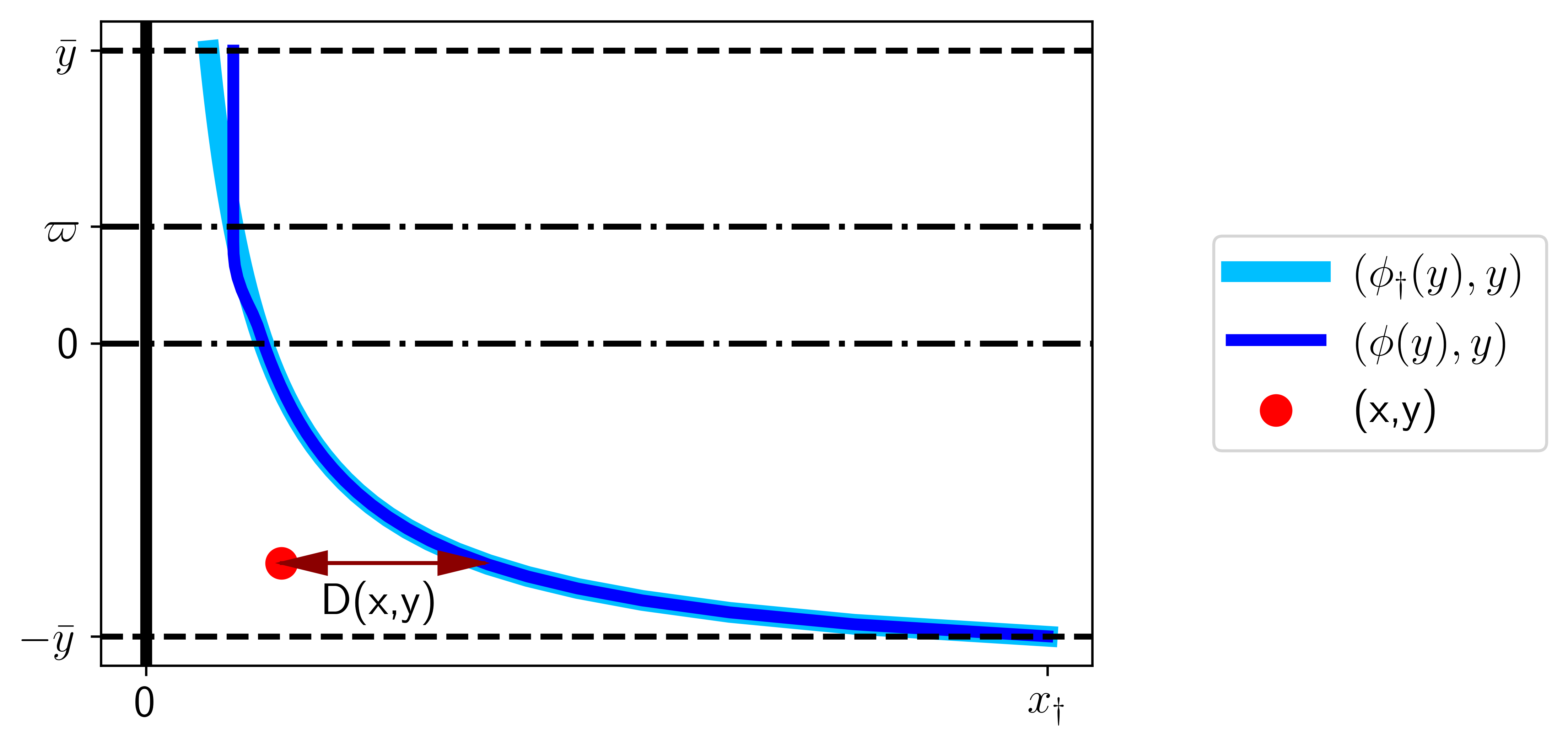}
    \caption{Danger function $D$}
    \label{fig:danger}
\end{figure}

Let's now define the function
\begin{equation*} D(x,y)\Def \phi(y)-x \end{equation*} 
for $(x,y)\in \R\times \R$, which measures the distance of the point $(x,y)$ to the left of the graph of $\phi$ and thus quantifies the \emph{danger} of a collision; see Figure \ref{fig:danger}.
\begin{lemma}\label{L:Danger} If $t<\tau^\eps$, $t\le \tau^{H,\eps}$, and $X^\eps_t\le \tfrac12 \uphi$, then $D(X^\eps_t,Y^\eps_t)\ge \tfrac12 \uphi$.
\end{lemma}
\begin{proof}
Under the stated assumptions, $|Y^\eps_t| \le \bar y$ and thus
\begin{equation*} D(X^\eps_t,Y^\eps_t)=\phi(Y^\eps_t)-X^\eps_t\ge \uphi - \tfrac12 \uphi = \tfrac 12\uphi. \end{equation*}
\end{proof}
\noindent This converts proximity to the boundary to values of $D$.  We will show that $D(X^\eps_t,Y^\eps_t)$ is likely remain small, implying that $X^\eps$ is unlikely to be small, i.e., collision is unlikely.

Let's next define
\begin{equation*} \Delta(x,y) \Def \left(\max\lb D(x,y),0\rb\right)^2 \end{equation*}
which focuses on the positive part of $D$.
We note that $\Delta$ is in $C^1(\R)$ and in fact has a piecewise-continuous second derivative.

Define now (similar to \eqref{E:stopTime_bound})
\begin{equation*} \tau^{D,\eps}=\inf\lb t\in [0,\tau^\eps): X^\eps_t<\tfrac12 \uphi\rb \qquad \text{(with $\inf\emptyset \Def \infty$)}. \end{equation*}
Similar to \eqref{E:aaa} and \eqref{E:aab},
\begin{equation*} \lb \tau^{D,\eps}\le t\rb = \bigcap_{\substack{t'>t \\ t'\in \Q}}\lb \bigcup_{\substack{s\in \Q \\ s<t'}}\bigcup_{m=1}^\infty \left(\lb X^\eps_s<\tfrac12 \uphi\rb \cap \{\tau^{\eps,\nicefrac{1}{m}}>s\}\right)\rb \in \filt_{t+}=\filt_t \end{equation*}
so $\tau^{D,\eps}$ is indeed a stopping time with respect to $\{\filt_t\}_{t\ge 0}$.

We can now prove
\begin{proposition}\label{P:longcollision}
We have that
\begin{equation*}
    \lim_{\substack{ \eps \to 0 \\ \eps \sqrt{L} \to 0}} \BP \lb \tau^{D,\eps} <L\wedge \tau^{H,\eps} \rb = 0.
\end{equation*}
\end{proposition}
\noindent Define
\begin{equation} \label{E:bardeltadef} \bar \delta_{\eqref{E:bardeltadef}}\Def \min\lb \frac{1}{2\bar y}, \frac{\uphi}{4}\rb.  \end{equation}
Then
\begin{equation*} \lb (x,y)\in \Space: \text{$x\ge \tfrac12 \uphi$ and $H(x,y)\le h_\circ+1$}\rb \subset \left[\tfrac12 \uphi,\infty\right)\times [-\bar y,\bar y]\subset \agree_{\bar \delta_{\eqref{E:bardeltadef}}}. \end{equation*}
Using the last equality of \eqref{E:LimitClaims},
\begin{equation}\label{E:finaleq} \tau^{D,\eps}\wedge \tau^{H,\eps}\le \tau^{\eps,\bar\delta_{\eqref{E:bardeltadef}}}. \end{equation} \begin{proof}[Proof of Proposition \ref{P:longcollision}]
By \eqref{E:finaleq}, we have that $\tau^{D,\eps}\wedge \tau^{H,\eps}\le \tau^{\eps,\bar\delta_{\eqref{E:bardeltadef}}}$ and we can use the second equality of \eqref{E:LimitClaims} setting $\delta$ to be $\bar \delta_{\eqref{E:bardeltadef}}$.

If $\tau^{D,\eps}<L\wedge \tau^{H,\eps}$, then $\tau^{D,\eps}<\infty$ and thus $X^\eps_{\tau^{D,\eps}\wedge \tau^{H,\eps}\wedge L}=\tfrac12 \uphi$, and consequently (by Lemma \ref{L:Danger})
\begin{equation*} D\left(X^\eps_{\tau^{D,\eps}\wedge \tau^{H,\eps}\wedge L},Y^\eps_{\tau^{D,\eps}\wedge \tau^{H,\eps}\wedge L}\right)\ge \tfrac12 \uphi. \end{equation*}
Thus
\begin{equation}\label{E:collection} \begin{aligned} \BP\lb \tau^{D,\eps}<L\wedge \tau^{H,\eps}\rb 
&\le \BP\lb \Delta\left(X^\eps_{\tau^{D,\eps}\wedge \tau^{H,\eps}\wedge L},Y^\eps_{\tau^{D,\eps}\wedge \tau^{H,\eps}\wedge L}\right)>\tfrac14 \uphi^2\rb \\
&\le \frac{4}{\uphi^2}\BE\left[\Delta\left(X^\eps_{\tau^{D,\eps}\wedge \tau^{H,\eps}\wedge L},Y^\eps_{\tau^{D,\eps}\wedge \tau^{H,\eps}\wedge L}\right)\right]. \end{aligned}\end{equation}

Let's now use Ito's formula on $\Delta$.  Since $x\mapsto (x^+)^2$ is not twice differentiable (thus precluding a direct application of Ito's rule), let's approximate.  Define
\begin{equation*} \smooth(x)\Def \begin{cases} 0 &\text{if $x<-1$} \\
2(x+1) &\text{if $-1\le x< 0$} \\
2 &\text{if $x\ge 0$}\end{cases}\end{equation*}
For $\vsig>0$, define
\begin{align*} \smooth_\vsig(x)&\Def \smooth\left(\nicefrac{x}{\vsig}\right)\\
\smooth^{(-1)}_\vsig(x)&\Def \int_{x_1=-\infty}^x \smooth_\vsig(x_1)dx_1 = \begin{cases} 0 &\text{if $x<-\vsig$} \\
\int_{x_1=-\vsig}^x \smooth_\vsig(x_1)dx_1 &\text{if $x\ge -\vsig$} \end{cases} \\
\smooth^{(-2)}_\vsig(x)&\Def \int_{x_2=-\infty}^x \lb \int_{x_1=-\infty}^{x_2} \smooth_\vsig(x_1)dx_1\rb dx_2 = \begin{cases} 0 &\text{if $x<-\vsig$} \\
\int_{x_1=-\vsig}^x (x-x_1)\smooth_\vsig(x_1)dx_1 &\text{if $x\ge -\vsig$}\end{cases} \end{align*}
pointwise (for $x\in \R$).
Then $\dot \smooth_\vsig^{(-2)}\equiv \smooth^{(-1)}_\vsig$ and $\ddot \smooth^{(-2)}_\vsig\equiv \smooth_\vsig$.  For $\vsig\in (0,1)$, we also have that
\begin{equation} \label{E:sigbnds} \begin{gathered} 0\le \smooth_\vsig(x)\le 2\bOne_{[-1,\infty)}(x)\\
0\le \smooth^{(-1)}_\vsig(x)\le 2(x+1)^+ \\
0\le \smooth^{(-2)}_\vsig(x)\le \lb (x+1)^+\rb^2 \\
\lim_{\vsig\searrow 0}\smooth_\vsig(x)= 2\bOne_{[0,\infty)}(x) \\
\lim_{\vsig\searrow 0}\smooth^{(-1)}_\vsig(x) = 2x^+ \\
\lim_{\vsig\searrow 0}\smooth^{(-2)}_\vsig(x) = \{x^+\}^2\end{gathered}\end{equation}
for each $x\in \R$.

Let's now use Ito's formula to write
\begin{align*}\label{E:myito_A}\smooth^{(-2)}_\vsig\left(\Delta\left(X^\eps_{\tau^{D,\eps}\wedge \tau^{H,\eps}\wedge L},Y^\eps_{\tau^{D,\eps}\wedge \tau^{H,\eps}\wedge L}\right)\right)&= \smooth^{(-2)}_\vsig\left(\Delta(x_\circ,y_\circ)\right)\\
&\qquad + \int_{s=0}^{\tau^{D,\eps}\wedge \tau^{H,\eps}\wedge L} \smooth^{(-1)}_\vsig\left(D(X^\eps_s,Y^\eps_s)\right)\digamma(X^\eps_s,Y^\eps_s)ds\\
&\qquad -\eps\int_{s=0}^{\tau^{D,\eps}\wedge \tau^{H,\eps}\wedge L} \smooth^{(-1)}_\vsig\left(D(X^\eps_s,Y^\eps_s)\right) \phi^{\prime}(Y^\eps_s)dW_s\\
&\qquad +\eps^2\int_{s=0}^{\tau^{D,\eps}\wedge \tau^{H,\eps}\wedge L} \smooth^{(-1)}_\vsig\left(D(X^\eps_s,Y^\eps_s)\right) \phi^{\prime \prime}(Y^\eps_s)ds\\
&\qquad +\tfrac12\eps^2\int_{s=0}^{\tau^{D,\eps}\wedge \tau^{H,\eps}\wedge L} \smooth_\vsig\left(D(X^\eps_s,Y^\eps_s)\right) \left(\phi^{\prime }(Y^\eps_s)\right)^2 ds
\end{align*}
where 
\begin{equation*} \digamma(x,y)\Def \phi'(y)\lb -\alpha\lb V\left(\frac{x}{d}\right)-v_\circ+y\rb -  \beta\frac{y}{x^2}\rb-y = -\phi'(y)\lb \alpha \lb V\left(\frac{x}{d}\right)-v_\circ\rb + \lb \alpha+\frac{\beta}{x^2}\rb y\rb -y\end{equation*}
for all $x>0$ and $y\in \R$.  Letting $\vsig\searrow 0$ and using \eqref{E:sigbnds}, the bounds of \eqref{E:phibds}, and dominated convergence, we get that
\begin{equation}\label{E:myito_B}\begin{aligned}\Delta\left(X^\eps_{\tau^{D,\eps}\wedge \tau^{H,\eps}\wedge L},Y^\eps_{\tau^{D,\eps}\wedge \tau^{H,\eps}\wedge L}\right)&= \Delta(x_\circ,y_\circ)\\
&\qquad + 2\int_{s=0}^{\tau^{D,\eps}\wedge \tau^{H,\eps}\wedge L}D(X^\eps_s,Y^\eps_s)\bOne_{[0,\infty)}(D(X^\eps_s,Y^\eps_ss)) \digamma(X^\eps_s,Y^\eps_s)ds\\
&\qquad - 2\eps\int_{s=0}^{\tau^{D,\eps}\wedge \tau^{H,\eps}\wedge L}D(X^\eps_s,Y^\eps_s)\bOne_{[0,\infty)}(D(X^\eps_s,Y^\eps_s)) \phi^{\prime}(Y^\eps_s)dW_s\\
&\qquad +\eps^2\int_{s=0}^{\tau^{D,\eps}\wedge \tau^{H,\eps}\wedge L}D(X^\eps_s,Y^\eps_s)\bOne_{[0,\infty)}(D(X^\eps_s,Y^\eps_s)) \phi^{\prime \prime}(Y^\eps_s)ds\\
&\qquad +\tfrac12\eps^2\int_{s=0}^{\tau^{D,\eps}\wedge \tau^{H,\eps}\wedge L}\bOne_{[0,\infty)}(D(X^\eps_s,Y^\eps_s)) \left(\phi^{\prime }(Y^\eps_s)\right)^2 ds
\end{aligned}\end{equation}

Optional sampling implies that
\begin{equation*} \BE\left[\int_{s=0}^{\tau^{D,\eps}\wedge \tau^{H,\eps}\wedge L}D(X^\eps_s,Y^\eps_s)\bOne_{[0,\infty)}(D(X^\eps_s,Y^\eps_s)) \phi^{\prime}(Y^\eps_s)dW_s\right]=0. \end{equation*}
Using \eqref{E:phibds}, we have that
\begin{equation*} \BE\left[\eps^2\int_{s=0}^{\tau^{D,\eps}\wedge \tau^{H,\eps}\wedge L}\bOne_{[0,\infty)}(D(X^\eps_s,Y^\eps_s)) \left(\phi^{\prime }(Y^\eps_s)\right)^2ds\right]\le \KK^2_{\eqref{E:phibds}}\eps^2 L. \end{equation*}
We also have that
\begin{equation*} D(X^\eps_s,Y^\eps_s)\le \phi(Y^\eps_s)\le \phi(-\bar y)=x_\dagger \end{equation*}
for $s\in [0,\tau^{D,\eps}\wedge \tau^{H,\eps}\wedge L)$ (since $\phi$ is decreasing on $[-\bar y,\bar y]$ by Lemma \ref{L:phiprops}, $X^\eps_s>0$ for $s\in [0,\tau^{D,\eps}\wedge \tau^{H,\eps}\wedge L)$, and using \eqref{E:stochasticbarrier}).  Thus
\begin{equation*} \BE\left[\left|\eps^2\int_{s=0}^{\tau^{D,\eps}\wedge \tau^{H,\eps}\wedge L}D(X^\eps_s,Y^\eps_s)\bOne_{[0,\infty)}(D(X^\eps_s,Y^\eps_s)) \phi^{\prime \prime}(Y^\eps_s)ds\right|\right]
\le x_\dagger\KK_{\eqref{E:phibds}}\eps^2 L. \end{equation*}

Let's finally look at the first term on the right of \eqref{E:myito_B}.  Fix $s\in [0,\tau^{D,\eps}\wedge \tau^{H,\eps}\wedge L)$ and assume that $D(X^\eps_s,Y^\eps_s)\ge 0$ (keeping in mind the term $\bOne_{[0,\infty)}(D(X^\eps_s,Y^\eps_s))$). We claim that then 
\begin{equation}\label{E:digammaclaim} \digamma(X^\eps_s,Y^\eps_s)\le 0. \end{equation}

Firstly, if $Y^\eps_s\ge \varpi_\dagger$, then $\phi'(Y^\eps_s)=0$ (recall \eqref{E:phibarder}) and $Y^\eps_s$ is positive, so
\begin{equation*}\digamma(X^\eps_s,Y^\eps_s)=-Y^\eps_s\le 0. \end{equation*}

We next observe that
\begin{equation}\label{E:hiho} X^\eps_s\le \phi(Y^\eps_s)\le \phi(-\bar y)=x_\dagger\le dV^{-1}\left(\nicefrac{v_\circ}{2}\right) \end{equation}
(since $D(X^\eps_s,Y^\eps_s)\ge 0$, $\phi$ is nonincreasing, using the initial condition \eqref{E:stochasticbarrier} of $\phi$, and then using the definition \eqref{E:xdaggerdef} of $x_\dagger$), thus implying that
\begin{equation}\label{E:Vsign} V\left(\frac{X^\eps_s}{d}\right)-v_\circ \le \tfrac12 v_\circ-v_\circ=-\tfrac12 v_\circ\le 0. \end{equation}

Assume next that $Y^\eps_s\le 0$.  From the first inequality of \eqref{E:hiho}, we have that
\begin{equation*} \alpha+\frac{\beta}{(X^\eps_s)^2}\ge \alpha+\frac{\beta}{\phi^2(Y^\eps_s)}. \end{equation*}
Using this and the last inequality of \eqref{E:Vsign} (and noting that $-\alpha \phi'(Y^\eps_s)\ge 0$ and $-Y^\eps_s\phi'(Y^\eps_s)\le 0$), we have that
\begin{align*} \digamma(X^\eps_s,Y^\eps_s)
&= -\alpha \phi'(Y^\eps_s)\lb V\left(\frac{X^\eps_s}{d}\right)-v_\circ \rb  -Y^\eps \lb \phi'(Y^\eps_s)\lb \alpha + \frac{\beta}{\left(X^\eps_s\right)^2}\rb + 1\rb \\
&\le 0-Y^\eps_s\lb \phi'(Y^\eps_s)\lb \alpha + \frac{\beta}{\phi^2(Y^\eps_s)}\rb + 1\rb =0\end{align*}
(using the fact that $\phi$ is nonincreasing and $Y^\eps_s$ is assumed to be negative), where we have used the ODE for $\phi_\dagger$ (i.e., the first case of \eqref{E:phibarder}). 

Thirdly, assume that $0\le Y^\eps_s\le \varpi_\dagger$.
In this case,
\begin{equation*} \alpha + \frac{\beta}{(X^\eps_s)^2}\le \alpha + \frac{\beta}{\nicefrac{\uphi^2}{4}} \end{equation*}
(Since $s\le \tau^{D,\eps}$, $X^\eps_s\ge \tfrac12 \uphi$), and thus
\begin{equation*} \digamma(X^\eps_s,Y^\eps_s)
\le -\phi'(Y^\eps_s)\lb \alpha\lb -\tfrac12 v_\circ\rb + \lb \alpha+ \frac{4\beta}{\uphi^2}\rb \varpi_\dagger \rb =0 \end{equation*}
(using the assumption that $Y^\eps_s$ is positive, the fact that $-\phi'(Y^\eps_s)>0$, using the first inequality of \eqref{E:Vsign}, and using the definition \eqref{E:varpidef} of $\varpi_\dagger$),

In light of \eqref{E:digammaclaim}, the claim follows from \eqref{E:collection}.
\end{proof}
We can now prove our main result.
\begin{proof}[Proof of Theorem \ref{T:main}]
Let's first write
\begin{equation}\label{E:finalineq} \BP\lb \tau^\eps<L\rb \le \BP\lb \tau^{H,\eps}\le L\rb + \BP\lb \tau^\eps<L,\, \tau^{H,\eps}>L\rb
\le \BP\lb \tau^{H,\eps}\le L\rb + \BP\lb \tau^\eps<L\wedge \tau^{H,\eps}\rb. \end{equation}

Let's bound the final term on the right of \eqref{E:finalineq}. Assume that $\tau^\eps<L\wedge \tau^{H,\eps}$.  From \eqref{E:finaleq}, the last equality of \eqref{E:LimitClaims} and this assumption, we have the chain
\begin{equation}\label{E:bound_prob_tau}
    \tau^{H,\eps} \wedge \tau^{D,\eps} < \tau^{\eps,\bar \delta_{\eqref{E:bardeltadef}}} < \tau^\eps <\tau^{H,\eps} \wedge L
\end{equation}
of inequalities.

If $\tau^{D,\eps}\ge \tau^{H,\eps}$, then $\tau^{H,\eps}=\tau^{H,\eps}\wedge \tau^{D,\eps}<\tau^{H,\eps}\wedge L$, which is contradiction.  Thus $\tau^{D,\eps}<\tau^{H,\eps}$ and hence
\begin{equation*}\tau^{D,\eps}=\tau^{H,\eps}\wedge \tau^{D,\eps}<\tau^{H,\eps}\wedge L
\end{equation*}

Combining things together, we have that
\begin{align*} \BP\lb \tau^\eps<L\rb &\le \BP\lb \tau^{H,\eps}\le L\rb +  \BP\lb \tau^{D,\eps}<L\wedge \tau^{H,\eps}\rb. \end{align*}
Using Theorem \ref{T:Hbounded} and Proposition \ref{P:longcollision}, the result follows.
\end{proof}

\section{Extensions}
Our analysis of the simplified and stylized noisy leader-follower dynamics of \eqref{E:initial_stochastic} can hopefully be extended to shed light on a number of more complex problems.  Dynamics near collisions for platoons might be considered by first of all adding some small time-varying noise to the lead velocity $v_\circ$, and then using information about how noise propagates through a pair of vehicles to a an entire platoon.  Mathematically some ideas from the theory of time-varying Hamiltonians might be used in place of \eqref{E:energyLevel}. Some ideas from string stability calculations (see the references listed in Section \ref{S:ModelsAndContributions}) might organize these calculations.  Most of our deterministic and stochastic boundary-analyses are likely to have a higher-dimensional counterpart; only the Poincar\'e-Bendixson calculations of Theorem \ref{T:absorbing} is likely to require a major overhaul (and that part of Theorem  \ref{T:absorbing} is in fact a statement of long-term stability).

One might also develop an abstract analysis of collision-avoidance terms to understand broad classes of collision-avoidance dynamics.  The work of  \cite{aghabayk2015state, gazis1961nonlinear} might serve as a guide in such an undertaking.

Finally, we believe that our main stochastic result, Theorem \ref{T:main}, might be tightened.  Our deterministic analysis of Section \ref{S:deterministic_collision}, and in particular Theorem \ref{T:absorbing}, preclude collision in finite time.  Continuity with respect to the $\eps$ of noise in \eqref{E:initial_stochastic} thus implies bounds on the likelihood of collision in any finite time interval, and large-deviations type estimates should be available.  The presence of a stable point (i.e., the final claim of Theorem \ref{T:absorbing}) suggests that collision times should be exponentially large \cite{MR1652127}.  Part of the motivation of our analysis, by comparison with a large-deviations approach, is a collection of boundary-region calculations which can be used when the follower vehicle is close to the lead vehicle.  A more global analysis which would combine large-deviations type asymptotics and boundary analysis near collisions might be of independent interest.  


\bibliographystyle{amsalpha}
\bibliography{reference}

\newcommand{\etalchar}[1]{$^{#1}$}
\providecommand{\bysame}{\leavevmode\hbox to3em{\hrulefill}\thinspace}
\providecommand{\MR}{\relax\ifhmode\unskip\space\fi MR }
\providecommand{\MRhref}[2]{%
  \href{http://www.ams.org/mathscinet-getitem?mr=#1}{#2}
}
\providecommand{\href}[2]{#2}
\begin{thebibliography}{SCDM{\etalchar{+}}18}

\bibitem[ASY15]{aghabayk2015state}
Kayvan Aghabayk, Majid Sarvi, and William Young, \emph{A state-of-the-art
  review of car-following models with particular considerations of heavy
  vehicles}, Transport reviews \textbf{35} (2015), no.~1, 82--105.

\bibitem[BHN{\etalchar{+}}94]{bando1994structure}
M~Bando, K~Hasebe, A~Nakayama, A~Shibata, and Y~Sugiyama, \emph{Structure
  stability of congestion in traffic dynamics}, Japan Journal of Industrial and
  Applied Mathematics \textbf{11} (1994), no.~2, 203.

\bibitem[BHN{\etalchar{+}}95]{bando1995dynamical}
Masako Bando, Katsuya Hasebe, Akihiro Nakayama, Akihiro Shibata, and Yuki
  Sugiyama, \emph{Dynamical model of traffic congestion and numerical
  simulation}, Physical review E \textbf{51} (1995), no.~2, 1035.

\bibitem[BMSN18]{bischoff2018autonomous}
Joschka Bischoff, Michal Maciejewski, Tilmann Schlenther, and Kai Nagel,
  \emph{Autonomous vehicles and their impact on parking search}, IEEE
  Intelligent Transportation Systems Magazine \textbf{11} (2018), no.~4,
  19--27.

\bibitem[CGYW19]{chen2019quantifying}
Yuche Chen, Jeffrey Gonder, Stanley Young, and Eric Wood, \emph{Quantifying
  autonomous vehicles national fuel consumption impacts: A data-rich approach},
  Transportation Research Part A: Policy and Practice \textbf{122} (2019),
  134--145.

\bibitem[CHM58]{chandler1958traffic}
Robert~E Chandler, Robert Herman, and Elliott~W Montroll, \emph{Traffic
  dynamics: studies in car following}, Operations research \textbf{6} (1958),
  no.~2, 165--184.

\bibitem[CL55]{MR0069338}
Earl~A. Coddington and Norman Levinson, \emph{Theory of ordinary differential
  equations}, McGraw-Hill Book Co., Inc., New York-Toronto-London, 1955.
  \MR{0069338}

\bibitem[CPT21]{chiarello2021statistical}
Felisia~Angela Chiarello, Benedetto Piccoli, and Andrea Tosin, \emph{A
  statistical mechanics approach to macroscopic limits of car-following traffic
  dynamics}, International Journal of Non-Linear Mechanics \textbf{137} (2021),
  103806.

\bibitem[CSSW17]{cui2017stabilizing}
Shumo Cui, Benjamin Seibold, Raphael Stern, and Daniel~B Work,
  \emph{Stabilizing traffic flow via a single autonomous vehicle: Possibilities
  and limitations}, Intelligent Vehicles Symposium (IV), 2017 IEEE, IEEE, 2017,
  pp.~1336--1341.

\bibitem[Dav03]{davis2003modifications}
LC~Davis, \emph{Modifications of the optimal velocity traffic model to include
  delay due to driver reaction time}, Physica A: Statistical Mechanics and its
  Applications \textbf{319} (2003), 557--567.

\bibitem[Dav14]{davis2014nonlinear}
\bysame, \emph{Nonlinear dynamics of autonomous vehicles with limits on
  acceleration}, Physica A: Statistical Mechanics and its Applications
  \textbf{405} (2014), 128--139.

\bibitem[dJF96]{de1996theory}
Eduardus~Marie de~Jager and JF~Furu, \emph{The theory of singular
  perturbations}, Elsevier, 1996.

\bibitem[EK09]{ethier2009markov}
Stewart~N Ethier and Thomas~G Kurtz, \emph{Markov processes: characterization
  and convergence}, John Wiley \& Sons, 2009.

\bibitem[FW98]{MR1652127}
M.~I. Freidlin and A.~D. Wentzell, \emph{Random perturbations of dynamical
  systems}, second ed., Grundlehren der Mathematischen Wissenschaften
  [Fundamental Principles of Mathematical Sciences], vol. 260, Springer-Verlag,
  New York, 1998, Translated from the 1979 Russian original by Joseph Sz\"ucs.
  \MR{MR1652127 (99h:60128)}

\bibitem[GGLP20]{garavello2020multiscale}
Mauro Garavello, Paola Goatin, Thibault Liard, and Benedetto Piccoli, \emph{A
  multiscale model for traffic regulation via autonomous vehicles}, Journal of
  Differential Equations \textbf{269} (2020), no.~7, 6088--6124.

\bibitem[GGS{\etalchar{+}}20]{gunter2020commercially}
George Gunter, Derek Gloudemans, Raphael~E Stern, Sean McQuade, Rahul Bhadani,
  Matt Bunting, Maria~Laura Delle~Monache, Roman Lysecky, Benjamin Seibold,
  Jonathan Sprinkle, et~al., \emph{Are commercially implemented adaptive cruise
  control systems string stable?}, IEEE Transactions on Intelligent
  Transportation Systems \textbf{22} (2020), no.~11, 6992--7003.

\bibitem[GHR61]{gazis1961nonlinear}
Denos~C Gazis, Robert Herman, and Richard~W Rothery, \emph{Nonlinear
  follow-the-leader models of traffic flow}, Operations research \textbf{9}
  (1961), no.~4, 545--567.

\bibitem[GP06]{garavello2006traffic}
Mauro Garavello and Benedetto Piccoli, \emph{Traffic flow on networks}, vol.~1,
  American institute of mathematical sciences Springfield, 2006.

\bibitem[GS15]{greenblatt2015automated}
Jeffery~B Greenblatt and Susan Shaheen, \emph{Automated vehicles, on-demand
  mobility, and environmental impacts}, Current sustainable/renewable energy
  reports \textbf{2} (2015), no.~3, 74--81.

\bibitem[HM08]{hamdar2008existing}
Samer~H Hamdar and Hani~S Mahmassani, \emph{From existing accident-free
  car-following models to colliding vehicles: exploration and assessment},
  Transportation research record \textbf{2088} (2008), no.~1, 45--56.

\bibitem[HT98]{helbing1998generalized}
Dirk Helbing and Benno Tilch, \emph{Generalized force model of traffic
  dynamics}, Physical review E \textbf{58} (1998), no.~1, 133.

\bibitem[JN03]{jost2003probabilistic}
Dominic Jost and Kai Nagel, \emph{Probabilistic traffic flow breakdown in
  stochastic car-following models}, Transportation Research Record
  \textbf{1852} (2003), no.~1, 152--158.

\bibitem[JWZ01]{jiang2001full}
Rui Jiang, Qingsong Wu, and Zuojin Zhu, \emph{Full velocity difference model
  for a car-following theory}, Physical Review E \textbf{64} (2001), no.~1,
  017101.

\bibitem[KC11]{kavathekar2011vehicle}
Pooja Kavathekar and YangQuan Chen, \emph{Vehicle platooning: A brief survey
  and categorization}, International Design Engineering Technical Conferences
  and Computers and Information in Engineering Conference, vol. 54808, 2011,
  pp.~829--845.

\bibitem[Kra98]{krauss1998microscopic}
Stefan Krau{\ss}, \emph{Microscopic modeling of traffic flow: Investigation of
  collision free vehicle dynamics}.

\bibitem[KS58]{kometani1958stability}
EIJI Kometani and Tsuna Sasaki, \emph{On the stability of traffic flow
  (report-i)}, Journal of the Operations Research Society of Japan \textbf{2}
  (1958), no.~1, 11--26.

\bibitem[LHZ21]{le2021air}
Zoe Le~Hong and Naomi Zimmerman, \emph{Air quality and greenhouse gas
  implications of autonomous vehicles in vancouver, canada}, Transportation
  Research Part D: Transport and Environment \textbf{90} (2021), 102676.

\bibitem[LLH17]{lu2017platooning}
Duo Lu, Zhichao Li, and Dijiang Huang, \emph{Platooning as a service of
  autonomous vehicles}, 2017 IEEE 18th International Symposium on A World of
  Wireless, Mobile and Multimedia Networks \textup{(}WoWMoM\textup{)}, IEEE,
  2017, pp.~1--6.

\bibitem[LTZ14]{laval2014parsimonious}
Jorge~A Laval, Christopher~S Toth, and Yi~Zhou, \emph{A parsimonious model for
  the formation of oscillations in car-following models}, Transportation
  Research Part B: Methodological \textbf{70} (2014), 228--238.

\bibitem[Mat21]{MatinThesis}
Hossein Nick~Zinat Matin, \emph{Asymptotic behavior of stochastic optimal
  velocity dynamical model}, Ph.D. thesis, University of Illinois, 2021.

\bibitem[MKL05]{mahnke2005probabilistic}
Reinhard Mahnke, Jevgenijs Kaupu{\v{z}}s, and I~Lubashevsky,
  \emph{Probabilistic description of traffic flow}, Physics Reports
  \textbf{408} (2005), no.~1-2, 1--130.

\bibitem[MM07]{mitra2007pollution}
Debojyoti Mitra and Asis Mazumdar, \emph{Pollution control by reduction of drag
  on cars and buses through platooning}, International Journal of Environment
  and Pollution \textbf{30} (2007), no.~1, 90--96.

\bibitem[Mor88]{morita1988brownian}
Akio Morita, \emph{Brownian motion of two interacting particles through
  logarithmic and lennard-jones potentials}, Journal of the Chemical Society,
  Faraday Transactions 2: Molecular and Chemical Physics \textbf{84} (1988),
  no.~11, 1855--1866.

\bibitem[MS20a]{matin2020nonlineara}
Hossein Nick~Zinat Matin and Richard~B Sowers, \emph{Nonlinear optimal velocity
  car following dynamics (i): Approximation in presence of deterministic and
  stochastic perturbations}, 2020 American Control Conference (ACC), IEEE,
  2020, pp.~410--415.

\bibitem[MS20b]{matin2020nonlinearb}
\bysame, \emph{Nonlinear optimal velocity car following dynamics (ii): Rate of
  convergence in the presence of fast perturbation}, 2020 American Control
  Conference (ACC), IEEE, 2020, pp.~416--421.

\bibitem[NLT{\etalchar{+}}19]{ngoduy2019langevin}
D~Ngoduy, S~Lee, M~Treiber, M~Keyvan-Ekbatani, and HL~Vu, \emph{Langevin method
  for a continuous stochastic car-following model and its stability
  conditions}, Transportation Research Part C: Emerging Technologies
  \textbf{105} (2019), 599--610.

\bibitem[Oks13]{oksendal2013stochastic}
Bernt Oksendal, \emph{Stochastic differential equations: an introduction with
  applications}, Springer Science \& Business Media, 2013.

\bibitem[SCC{\etalchar{+}}19]{stern2019quantifying}
Raphael~E Stern, Yuche Chen, Miles Churchill, Fangyu Wu, Maria~Laura
  Delle~Monache, Benedetto Piccoli, Benjamin Seibold, Jonathan Sprinkle, and
  Daniel~B Work, \emph{Quantifying air quality benefits resulting from few
  autonomous vehicles stabilizing traffic}, Transportation Research Part D:
  Transport and Environment \textbf{67} (2019), 351--365.

\bibitem[SCDM{\etalchar{+}}18]{stern2018dissipation}
Raphael~E Stern, Shumo Cui, Maria~Laura Delle~Monache, Rahul Bhadani, Matt
  Bunting, Miles Churchill, Nathaniel Hamilton, R'mani Haulcy, Hannah Pohlmann,
  Fangyu Wu, Benedetto Piccoli, Benjamin Siebold, Jonathan Sprinkle, and
  Daniel~B. Work, \emph{Dissipation of stop-and-go waves via control of
  autonomous vehicles: Field experiments}, Transportation Research Part C:
  Emerging Technologies \textbf{89} (2018), 205--221.

\bibitem[SH07]{schonhof2007empirical}
Martin Sch{\"o}nhof and Dirk Helbing, \emph{Empirical features of congested
  traffic states and their implications for traffic modeling}, Transportation
  Science \textbf{41} (2007), no.~2, 135--166.

\bibitem[Sow05]{sowers2005boundary}
Richard~B Sowers, \emph{A boundary layer theory for diffusively perturbed
  transport around a heteroclinic cycle}, Communications on Pure and Applied
  Mathematics: A Journal Issued by the Courant Institute of Mathematical
  Sciences \textbf{58} (2005), no.~1, 30--84.

\bibitem[TK13]{treiber2013traffic}
Martin Treiber and Arne Kesting, \emph{Traffic flow dynamics}, Traffic Flow
  Dynamics: Data, Models and Simulation, Springer-Verlag Berlin Heidelberg
  (2013), 983--1000.

\bibitem[TKH10]{treiber2010three}
Martin Treiber, Arne Kesting, and Dirk Helbing, \emph{Three-phase traffic
  theory and two-phase models with a fundamental diagram in the light of
  empirical stylized facts}, Transportation Research Part B: Methodological
  \textbf{44} (2010), no.~8-9, 983--1000.

\bibitem[TM15]{talebpour2015influence}
Alireza Talebpour and Hani~S Mahmassani, \emph{Influence of autonomous and
  connected vehicles on stability of traffic flow}, Tech. report, 2015.

\bibitem[TS14]{tordeux2014collision}
Antoine Tordeux and Armin Seyfried, \emph{Collision-free nonuniform dynamics
  within continuous optimal velocity models}, Physical Review E \textbf{90}
  (2014), no.~4, 042812.

\bibitem[Wil99]{willis1999proportional}
MJ~Willis, \emph{Proportional-integral-derivative control}, Dept. of Chemical
  and Process Engineering University of Newcastle (1999).

\bibitem[WW11]{wilson2011car}
R~Eddie Wilson and Jonathan~A Ward, \emph{Car-following models: fifty years of
  linear stability analysis--a mathematical perspective}, Transportation
  Planning and Technology \textbf{34} (2011), no.~1, 3--18.

\bibitem[YLK{\etalchar{+}}18]{yuan2018geometric}
Kai Yuan, Jorge Laval, Victor~L Knoop, Rui Jiang, and Serge~P Hoogendoorn,
  \emph{A geometric brownian motion car-following model: towards a better
  understanding of capacity drop}, Transportmetrica B: Transport Dynamics
  (2018).

\bibitem[YY19]{ye2019evaluating}
Lanhang Ye and Toshiyuki Yamamoto, \emph{Evaluating the impact of connected and
  autonomous vehicles on traffic safety}, Physica A: Statistical Mechanics and
  its Applications \textbf{526} (2019), 121009.

\bibitem[ZSSB19]{zeng2019joint}
Tengchan Zeng, Omid Semiari, Walid Saad, and Mehdi Bennis, \emph{Joint
  communication and control for wireless autonomous vehicular platoon systems},
  IEEE Transactions on Communications \textbf{67} (2019), no.~11, 7907--7922.

\end{thebibliography}
\end{document}